\documentclass[a4paper]{article}

\usepackage{amsmath,amssymb,amsthm,amsfonts,amscd,euscript} 

\theoremstyle{plain}
\newtheorem{theorem}{Theorem}[section]
\newtheorem{propos}{Proposition}[section]
\newtheorem{corr}{Corollary}[section]
\theoremstyle{definition}

\theoremstyle{remark}
\newtheorem{note}{Note}[section]

\newcommand{\Hom}{\mathop{\rm Hom}\nolimits}

\newcommand{\Lk}{\mathop{\rm Lk}\nolimits}

\newcommand{\bideg}{\mathop{\rm bideg}\nolimits}

\newcommand{\den}{\mathop{\rm den}\nolimits}

\newcommand{\St}{\mathop{\rm Star}\nolimits}

\author{Alexander A. Gaifullin}

\title{Local formulae for combinatorial Pontrjagin classes}

\date{}
\begin{document}

\maketitle

\begin{abstract}
By $p(|K|)$ denote the characteristic class of a combinatorial
manifold $K$ given by the polynomial $p$ in Pontrjagin classes of
$K$. We prove that for any polynomial $p$ there exists a function
taking each combinatorial manifold $K$ to a  {\it rational simplicial 
cycle} $z(K)$ such that: (1) the
Poincar\'e dual of $z(K)$ represents the cohomology class
$p(|K|)$; (2) a coefficient of each simplex $\Delta$ in the cycle
$z(K)$ is determined only by the combinatorial type of $\Lk\Delta$.
We also prove that if a function $z$ satisfies the condition (2), then
this function automatically satisfies the condition (1) for some polinomial $p$.
We describe explicitly all such functions $z$ for the first
Pontrjagin class. We obtain estimates for denominators of 
coefficients of simplices in the cycles $z(K)$.
\end{abstract}

\section{Introduction}

The following general problem was studied by many researchers,
see~\cite{C,GGL1,GGL2,GMP,LR}. Given a triangulation of a
manifold one need to construct a cycle, whose Poincar\'e dual
represents a given Pontrjagin class of this manifold. In
addition, one usually wants the coefficient of each simplex
in this cycle to be determined by the structure of the manifold in
some neighborhood of the simplex. First let us discuss the most
important results concerning this problem.

A. M. Gabrielov, I. M. Gelfand and M. V. Losik~\cite{GGL1,GGL2}
found an explicit formula for the first rational Pontrjagin class
of a smooth manifold. To apply this formula one need the
manifold to be endowed with a smooth triangulation satisfying some
special condition. In their paper~\cite{GMP}, I. M. Gelfand and
R. D. MacPherson considered simplicial manifolds endowed with some
additional combinatorial data called {\it a fixing cycle}. A
fixing cycle is a combinatorial analogue of a smooth structure
and can be induced by a given smooth structure. For simplicial
manifolds with a given fixing cycle I. M. Gelfand and R. D.
MacPherson constructed rational cycles whose Poincar\'e dual
represent the normal Pontrjagin classes of the manifolds. In
these cycles the coefficient of simplex depends both on the
combinatorial structure of the neighborhood of this simplex and on
the restriction of the fixing cycle to this neighborhood.

The other approach is Cheeger's. In~\cite{C} he obtained 
explicit formulae for cycles whose Poincar\'e dual represent
{\it real} Pontrjagin classes. In these cycles the coefficient of
simplex depends only on the combinatorial type of the link of
this simplex. Cheeger's formula can be applied for any
pseudomanifold. It is unknown if the obtained cycles are rational.

Assume that for any combinatorial manifold $K$ we have a {\it
cycle} in the cooriented simplicial chains of $K$ given by
$$
z(K)=\sum_{\Delta\in K,\,\dim\Delta=\dim K-n}f(\Lk\Delta)\Delta
$$
where the value $f(L)$ is determined by the isomorphism class of
the oriented $(n-1)$-dimensional PL sphere $L$ and the function
$f$ does not depend on the manifold $K$. Then we say that $z$ is
{\it a characteristic local cycle} of codimension
$n$. The function $f$ is called {\it a local formula} for this 
characteristic local cycle. We prove that for any rational characteristic class
$p\in H^*(BPL;\mathbb{Q})$ there exists a {\it rational} characteristic 
local cycle
$z$ such that the Poincar\'e dual of the cycle $z(K)$ represents
the cohomology class $p(|K|)$ for any combinatorial manifold $K$.
($BPL$ is the classifying space for stable PL bundles.) This
improves a theorem of N. Levitt and C. Rourke~\cite{LR}. They
obtained a similar result for the cycles given by
$$
z(K)=\sum_{\Delta\in K,\,\dim\Delta=\dim K-n}g(\Lk\Delta,\dim
K)\Delta
$$
Such cycles are not characteristic local cycles because the function $g$ depends
on the dimension of $K$.

In section 2 we define the cochain complex $\mathcal{T}^*(\mathbb
Q)$ whose elements are skew-symmetric functions from the set of all
isomorphism classes of oriented PL spheres to $\mathbb Q$. We prove
that the function $f\in\mathcal{T}^*(\mathbb Q)$ is a cocycle iff
it is a local formula for a characteristic local cycle. We prove that there
exists an isomorphism $H^*(\mathcal T^*(\mathbb Q))\cong
H^*(BPL;\mathbb Q)=H^*(BO;\mathbb Q)$. Hence any rational characteristic local
cycle represents the homology classes dual to some polynomial in
the Pontrjagin classes of the manifolds. Besides, a local formula
for a given rational characteristic class is unique up to a coboundary.

In section 3 we obtain an explicit formula for all rational characteristic local
cycles $z$ such that the Poincar\'e dual of $z(K)$ represents the
first Pontrjagin class of a combinatorial manifold $K$. This
result is new because the formulae of~\cite{GGL1,GGL2,GMP} cannot
be applied for an arbitrary combinatorial manifold and the
formulae of~\cite{C} give only {\it real} characteristic local cycles. We use
the following approach. First we find explicitly all rational characteristic
local cycles of codimension $4$, i. e. all local formulae
$f\in\mathcal T^4(\mathbb Q)$. Then we prove and use the theorem that any
such characteristic local cycle represents the homology classes dual to  the
first Pontrjagin class multiplied by some rational constant. The
usage of {\it bistellar moves} is very important.

In section 3.6 we find some estimates for the denominators of the
local formulae's values. In particular, we prove that if $f$ is
a local formula for the first Pontrjagin class and $q$ be an
integer number, then there exists a PL sphere $L$ such that the
denominator of $f(L)$ is divisible by $q$. Hence there exist no
integer local formulae representing nontrivial homology classes.

A brief formulation of our results is contained in~\cite{G1}. In
this paper we omit proofs of several propositions. Full proofs
can be found in~\cite{G2}.

{\bf Terminology and notation.} All necessary definitions and
results of PL topology can be found in~\cite{RS}. In the sequel
all manifolds, triangulations and bordisms are supposed to be
piecewise linear. All manifolds are supposed to be closed. All
bordisms are supposed to be oriented. An isomorphism of oriented
simplicial complexes is an isomorphism preserving orientation. An
isomorphism changing orientation is called an antiisomorphism.
Let $K$ be a simplicial complex on a set $S$. By $CK$, $\Sigma
K$, and $K'$ denote respectively the cone and the suspension over
$K$, and the barycentric subdivision of $K$. {\it The full
subcomplex} spanned by a subset $V\subset S$ is the subcomplex
$L\subset K$ consisting of all simplices $\Delta\in K$ such that
all vertices of $\Delta$ belong to $V$. By $\Lk\Delta$ and $\St\Delta$ denote 
respectively the link and the star of a simplex $\Delta$.

\section{Local formulae}

\subsection{Main definitions and results}

Let $\mathcal T_n$ be the set of all isomorphism classes of oriented
$(n-1)$-dimensional PL spheres. Usually we don't distinguish
between a PL sphere and its isomorphism class. For any $L\in\mathcal T_n$
by $-L$ denote the same triangulation with the opposite
orientation. Let $L\in\mathcal T_n$ be {\it symmetric} if there exists an
antiautomorphism of $L$. Let $G$ be an abelian group. By $\mathcal T^n(G)$
denote the abelian group of all functions $f:\mathcal T_n\to G$ such that
$f(L)=f(-L)$ for any $L\in\mathcal T_n$. We assume that $\mathcal T^0(G)=G$,
$\mathcal T^{-n}(G)=0$, $n>0$. Let the differential $\delta :\mathcal T^n(G)\to
\mathcal T^{n+1}(G)$ be given by $(\delta f)(L) = \sum f(\Lk v)$, where the
sum is over all vertices $v\in L$ and the orientation of $\Lk v$
is induced by the orientation of $L$. Evidently, $\delta^2=0$. Thus
$\mathcal T^*(G)$ is a cochain complex.

Let $K$ be an $m$-dimensional combinatorial manifold. Let $\widehat G$
be the local system on $|K|$ with fiber $G$ and twisting given by
the orientation. A coorientation of the simplex $\Delta^n\in K$
is an orientation of $\Lk\Delta^n$. Any $m$-simplex is supposed
to be positevely cooriented. By $\widehat C_n(K;G)$ denote the
complex of cooriented simplicial chains of $K$. Let $\widehat\partial
:\widehat C_n(K;G)\to \widehat C_{n-1}(K;G)$ be the boundary operator.
(The incidence coefficient of two simplices
$\tau^{k-1}\subset\sigma^k$ is equal to $+1$ iff the orientation
of $\Lk\sigma^k$ is induced by the orientation of
$\Lk\tau^{k-1}$.) The homology of $\widehat C_*(K;G)$ is equal to
$H_*(|K|;\widehat G)$. If $K$ is oriented then we have the
augmentation $\epsilon:\widehat C_0(K;G)\to G$.

Consider $f\in\mathcal T^n(G)$. By $f_{\sharp}(K)$ denote the cooriented
chain $f_{\sharp}(K) = \sum_{\Delta^{m-n}\in
K}f(\Lk\Delta^{m-n})\Delta^{m-n}\in\widehat C_{m-n}(K;G)$. Evidently
the addend $f(\Lk\Delta^{m-n})\Delta^{m-n}$ does not depend on the
coorientation of $\Delta^{m-n}$. Let $f\in\mathcal T^n(G)$ be {\it a
local formula} if for any combinatorial manifold $K$ the
cooriented chain $f_{\sharp}(K)$ is a cycle. {\sloppy

}

\begin{propos}
1) $f$ is a local formula if and only if $f$ is a cocycle in the
cochain complex $\mathcal T^*(G)$.\\
2) If $f$ is a coboundary in $\mathcal T^*(G)$, then for any
combinatorial manifold $K$ the cycle $f_{\sharp}(K)$ is a boundary.\\
3) Let $K_1$ and $K_2$ be two  triangulations of a manifold
$M^m$. If $f$ is a local formula, then $f_{\sharp}(K_1)$ and
$f_{\sharp}(K_2)$ are homologous.
\end{propos}

\begin{proof}
Notice that $\widehat\partial f_{\sharp}(K) = (\delta f)_{\sharp}(K)$
for any $f\in\mathcal T^n(G)$. The second claim of the proposition
follows. Also, if $f$ is a cocycle, then $f$ is a local formula.

Suppose $\delta f\ne 0$. Then there exists $L\in\mathcal T_{n+1}$ such
that $(\delta f)(L)\ne 0$. If $m>n$, then there exists an
$m$-dimensional combinatorial manifold $K$ such that
$\Lk\Delta\cong L$ for some simplex $\Delta\in K$. Then the
coefficient at the simplex $\Delta$ in the chain $\widehat\partial
f_{\sharp}(K) = (\delta f)_{\sharp}(K)$ is nonzero. Hence
$f_{\sharp}(K)$ is not a cycle. Therefore $f$ is not a local
formula.

Let us now prove the third claim of the proposition. Assume that
$m>n$. $K_1$ can be transformed to $K_2$ by the sequence of
stellar subdivisions and inverse stellar subdivisions
(see~\cite{A}). We can assume without loss of generality that
$K_1$ can be transformed to $K_2$ by a stellar subdivision of some
simplex $\Delta$. Then the support of the cycle
$f_{\sharp}(K_2)-f_{\sharp}(K_1)$ is contained in the subcomplex
$\St\Delta$, which is contractible. The proposition follows.

The case $n=m$ is proved in section 2.2.
\end{proof}

Consider $\psi\in H^n(\mathcal T^*(G))$. For any manifold $M^m$ by
$\psi_{\sharp}(M^m)\in H_{m-n}(M^m;\widehat G)$ denote the homology
class represented by $f_{\sharp}(K)$, where $K$ is an arbitrary
triangulation of $M^m$ and $f$ is an arbitrary
representative of $\psi$. By $\psi^{\sharp}(M^m)\in H^n(M^m;G)$
denote the Poincar\'e dual of $\psi_{\sharp}(M^m)$. If $M^m$ is oriented
and $m=n$, then let $\psi^{\star}(M^n)$ be given by $\psi^{\star}(M^n)=\langle
\psi^{\sharp}(M^n), [M^n]\rangle$.

\begin{theorem}
For any rational characteristic class $p\in H^n(BPL;\mathbb{Q})$
there exists a unique cohomology class $\phi_p\in
H^n(\mathcal{T}^*(\mathbb{Q}))$ such that
$\phi_p^{\sharp}(M)=p(M)$ for any manifold $M$. The homomorphism
$H^*(BPL;\mathbb{Q})\to H^*(\mathcal{T}^*(\mathbb{Q}))$ given by
$p\mapsto\phi_p$ is an isomorphism.
\end{theorem}

\begin{corr}
There is an additive isomorphism
$H^*(\mathcal{T}^*(\mathbb{Q}))\cong \mathbb{Q}[p_1,p_2,\ldots]$,
$\deg p_j=4j$.
\end{corr}

\begin{corr}
If $f\in\mathcal T^*(\mathbb Q)$ and $f_{\sharp}(K)$ 
is a boundary for any combinatorial manifold
$K$, then $f$ is a coboundary.
\end{corr}

\subsection{Invariance under bordisms}

\begin{propos}
Let $\psi\in H^n(\mathcal T^*(G))$ be an arbitrary cohomology class.
Then $\psi^{\star}(M_1^n)=\psi^{\star}(M_2^n)$ for any two bordant
oriented manifolds $M_1^n$ and $M_2^n$.
\end{propos}

Let $f\in\mathcal T^n(G)$ be a local formula. Let $L$ be an
$n$-dimensional null cobordant oriented combinatorial manifold.
Prove that $\epsilon(f_{\sharp}(L))=0$. Proposition 2.2 and the
unproved case of Proposition 2.1 follow.

Let $K$ be an $(n+1)$-dimensional combinatorial manifold with
boundary such that $\partial K=L$. Suppose $u$ is a cooriented
vertex of $L$; then $\Lk_Ku$ is a triangulation of $n$-disk. The
coorientation of the vertex $u$ in the simplicial complex $K$
induces the coorientation the vertex $u$ in the complex $L$.
Hence we have the monomorphism
$\widehat{C}_0(L;G)\to\widehat{C}_0(K;G)$. Obviously,
$\partial\Lk_Ku=\Lk_Lu$. By $\Lk^*_Ku$ denote the simplicial
complex $\Lk_Ku\cup_{\Lk_Lu}C(\Lk_Lu)$ whose orientation is
induced by the orientation of $\Lk_Ku$. Then 
$\Lk^*_Ku\in\mathcal T_{n+1}$. Similarly, for any cooriented 
1-simplex $e\in
L$ define the PL sphere $\Lk^*_Ke\in\mathcal T_n$.

Let the 1-chain $f_{\sharp}(K)\in\widehat C_1(K;G)$ be given by
$f_{\sharp}(K)=\sum f(\Lk^*_Ke)e+\sum f(\Lk_Ke)e$ where the first
sum is over all 1-simplices $e\in L$ and the second sum is over
all 1-simplices $e\in (K\backslash L)$. Since $f$ is a local
formula, we have $(\delta f)(\Lk^*_Ku)=0$ for any vertex $u\in L$.
Therefore $\sum f(\Lk_K^*e)+\sum f(\Lk_Ke)=f(\Lk_Lu)$, where the
first sum is over all 1-simplices $e\in L$ such that $u\in e$ and
the second sum is over all 1-simplices $e\in K\backslash L$ such
that $u\in e$ (the coorientation of $e$ is chosen so that the
incidence coefficient of the pair $(v,e)$ is equal to $+1$).
Consequently $\widehat\partial f_{\sharp}(K)=i(f_{\sharp}(L))$. Hence
$\epsilon(f_{\sharp}(L))=0$.

\begin{note}
Indeed, the formula $\widehat\partial
f_{\sharp}(K)=i(f_{\sharp}(\partial K))$ holds for a
combinatorial manifold $K$ of any dimension $m\ge n+1$.
\end{note}

\subsection{The isomorphism $\star$}
Let $\Omega_n$ be a group of oriented $n$-dimensional PL bordisms of a
point. From Proposition 2.2 it follows that the homomorphism
$\star:H^n(\mathcal T^*(G))\to \Hom(\Omega_n,G)$ taking each cohomology
class $\psi$ to $\psi^{\star}$ is well defined.  There exists a
canonical isomorphism $\Hom (\Omega_n,\mathbb{Q})\cong
H^n(BPL;\mathbb{Q})$. Hence we have a homomorphism $\sharp
:H^n(\mathcal{T}^*(\mathbb{Q}))\to H^n(BPL;\mathbb{Q})$.

\begin{theorem}
The homomorphism
$\star:H^n(\mathcal{T}^*(\mathbb{Q}))\to\Hom(\Omega_n,\mathbb{Q})$
is an isomorphism.
\end{theorem}

\begin{corr}
The homomorphism $\sharp:H^n(\mathcal{T}^*(\mathbb{Q}))\to
H^n(BPL;\mathbb{Q})$ is an isomorphism.
\end{corr}

The fact that $\star$ is an epimorphism can be easily deduced
from the results of Levitt and Rourke~\cite{LR}. Their approach
uses an explicit construction of the classifying space
$B\widetilde{PL}_n$ for block bundles. However it is possible to
prove that $\star$ is an isomorphism quite elementary,
see~\cite{G2}. Here we omit the proof of the fact that $\star$ is
a monomorphism, see~\cite{G2}.

\subsection{The cohomology classes $\psi^{\sharp}(M^m)$ for
manifolds of an arbitrary dimension}

From Proposition 2.2 it follows that for any cohomology class
$\psi\in H^n(\mathcal{T}^*(\mathbb Q))$ there exists a rational
characteristic class $p=\sharp(\psi)$ such that
$\psi^{\sharp}(M^n)=p(M^n)$ for any $n$-dimensional manifold
$M^n$. Let us prove that this formula holds for a manifold of any
dimension.

\begin{propos}
$\psi^{\sharp}(M^m)=p(M^m)=\sharp(\psi)(M^m)$ for any manifold
$M^m$, $m\ge n$.
\end{propos}

Obviously, Theorem 2.1 follows from Corollary 2.3 and Proposition
2.3.

Let $P$ be a compact polyhedron. Let $Q\subset P$ be a closed PL
subset. We say that $P$ is a manifold with singularities in $Q$
if $P\setminus Q$ is a (nonclosed) manifold. We have the
Lefschetz duality $H_{m-n}(P,Q;\widehat{G})\cong H^n(P\setminus
Q;G)$ for $n<m$. Let $\psi\in H^n(\mathcal{T}^*(\mathbb{Q}))$ be
a cohomology class satisfying $n<m$. Consider an arbitrary
triangulation $K$ of $P$. Let $L\subset K$ be the subcomplex
consisting of all closed simplices whose intersections with $Q$
are not empty. Let $f\in\mathcal{T}^n(G)$ be a local formula
representing $\psi$. Let the chain $f_{\sharp}(K,L)\in \widehat
C_{m-n}(K,L;G)$ be given by $\sum
f(\Lk\Delta^{m-n})\Delta^{m-n}$, where the sum is over all
$(m-n)$-simplices $\Delta^{m-n}\in K\setminus L$. Arguing as in
the proof of Proposition 2.1 it is easy to prove that
$f_{\sharp}(K,L)$ is a relative cycle, whose homology class does
not depend on the choice of a triangulation $K$ and a local
formula $f$ representing $\psi$. Thus the classes
$\psi_{\sharp}(P,Q)\in H_{m-n}(P,Q;\widehat{G})$ and
$\psi^{\sharp}(P,Q)\in H^n(P\backslash Q;G)$ are well defined.
Let $S\subset P\setminus Q$ be a compact subset. Evidently, the
cohomology class $\left.\psi^{\sharp}(P,Q)\right|_S$ is
determined by the topology of some neighborhood $U\supset S$.
Hence we have the following proposition.

\begin{propos}
Let $N^k$ be a (closed) oriented manifold. Let $\psi\in
H^n(\mathcal T^*(G))$, $n\le k$. Consider an $m$-dimensional
manifold $P$ with singularities in $Q$. Let $i:N^k\hookrightarrow
P\backslash Q$ be an embedding such that $i(N^k)\subset
P\setminus Q$ is a submanifold with trivial normal bundle.  Put
$\psi^{\sharp}_m(N^{k})= i^*(\psi^{\sharp}(K,L))$. Then the
cohomology class $\psi^{\sharp}_m(N^{k})$ depends only on the
manifold $N^k$ and the number $m$ and does not depend on the
choice of the triple $(P,Q,i)$.
\end{propos}

\begin{propos}
$\psi^{\sharp}_m(N^k)=\psi^{\sharp}(N^k)$ for any $m>k$.
\end{propos}

\begin{proof}
The join $N^k*\Delta^{m-k-1}$ is an $m$-dimensional manifold with
singularities in $N^k\sqcup\Delta^{m-k-1}$. Points of
$N^k*\Delta^{m-k-1}$ are linear combinations
$t_0x+\sum_{j=1}^{m-k}t_jy_j$, where $y_1,y_2,\ldots,y_{m-k}$ are
the vertices of the simplex $\Delta^{m-k-1}$, $x\in N^k$, $t_j\ge
0$, $j=0,1,\ldots,m-k$, $\sum_{j=0}^{m-k}t_j=1$. Let
$i:N^k\hookrightarrow N^k*\Delta^{m-k-1}$ be the embedding given
by $i(x)=\frac{1}{m-k+1}(x+\sum_{j=1}^{m-k}y_j)$. Then $i(N^k)$
is a submanifold with trivial normal bundle. Let $K$ be an
arbitrary triangulation of $N^k$. The embedding $i$ is
transversal to simplices of $K$. We have
$|\tau*\Delta^{m-k-1}|\cap i(N^k)=i(|\tau|)$ and
$\Lk_{K*\Delta^{m-k-1}}(\tau*\Delta^{m-k-1})=\Lk_K\tau$ for any
simplex $\tau\in K$. Hence for any local formula $f$ the
intersection of the cycles
$f_{\sharp}(K*\Delta^{m-k-1},K\sqcup\Delta^{m-k-1})$ and
$i_*([N^k])$ coincides with the cycle $i_*(f_{\sharp}(K))$.
Therefore
$\psi^{\sharp}_m(N^k)=\left.\psi^{\sharp}(M^m)\right|_{N^k}=\psi^{\sharp}(N^k)$.
\end{proof}

\begin{proof}[Proof of Proposition 2.3]
Consider the case $G=\mathbb Q$, $n=k$. Let $M^m$ be an oriented
manifold such that $m>n$. By Propositions 2.4 and 2.5 we have
$\left.\psi^{\sharp}(M^m)\right|_{N^n}=p(N^n)$ for any
submanifold $N^n\subset M^m$ with trivial normal bundle. From the
results of Rokhlin, Schwarz and Thom it follows that
$\psi^{\sharp}(M^m)=p(M^m)$ if $m>2n+1$.

Assume that $n<m\le 2n+1$. Let $i:M^m\hookrightarrow M^m\times
S^{n+1}$ be the standard embedding. By Propositions 2.4 and 2.5
we have $i^*(\psi^{\sharp}(M^m\times
S^{n+1}))=\psi^{\sharp}(M^m)$. On the other hand, $i^*(p(M^m\times
S^{n+1}))=p(M^m)$ and $\psi^{\sharp}(M^m\times
S^{n+1})=p(M^m\times S^{n+1})$. Therefore
$\psi^{\sharp}(M^m)=p(M^m)$.

Assume now that $M^m$  is not orientable. Let $\pi
:\widetilde{M}^m\to M^m$ be the oriented two-fold covering. Then
$\pi^*(\psi^{\sharp}(M^m))=\psi^{\sharp}(\widetilde M^m)$,
$\pi^*(p(M^m))=p(\widetilde M^m)$. The proposition follows because
$\pi^*$ is a monomorphism.
\end{proof}

\section{Explicit local formulae for the first Pontrjagin classes}
\subsection{Graphs $\Gamma_n$ and their cohomology}

Let $K$ be a combinatorial manifold. Assume that there is a
simplex $\Delta_1\in K$ such that $\Lk\Delta_1=\partial\Delta_2$.
(Obviously, $\Delta_2\notin K$.) Then $\Delta_1*\partial\Delta_2$
is a full subcomplex of $K$. The {\it bistellar move} associating
with the simplex $\Delta_1$ is the operation taking $K$ to the
simplicial complex $(K\backslash (\Delta_1*\partial\Delta_2)
)\cup(\partial\Delta_1*\Delta_2)$. If $\dim\Delta=0$, then
$\partial\Delta=\emptyset$. We assume that
$\Delta*\emptyset=\Delta$ for any simplex $\Delta$. Hence the
bistellar move associating with a maximal simplex of $K$ is just
the stellar subdivision of this simplex. The bistellar move
associating with a $0$-simplex is the inverse stellar
subdivision. By a theorem of Pachner~\cite{P} (see also~\cite{BP})
two PL manifolds are PL homeomorphic if and only if there exists a
sequence of bistellar moves transforming the first manifold to
the second one.

Consider two PL spheres $L_1,L_2\in\mathcal T_{n+1}$. Let
$\beta_1$ and $\beta_2$ be two bistellar moves transforming $L_1$
to $L_2$ such that $\beta_1$ and $\beta_2$ are associated with
simplices $\Delta_1$ and $\Delta_2$ respectively. We say that
$\beta_1$ and $\beta_2$ are {\it equivalent} if there exists an
automorphism of $L_1$ taking the simplex $\Delta_1$ to
$\Delta_2$. In the sequel we do not distinguish between a
bistellar move and its equivalence class. For any bistellar move
$\beta$ by $\beta^{-1}$ denote the inverse bistellar move. If a 
bistellar move $\beta$ transforming some PL sphere $L$ to itself
is equivalent to the bistellar move $\beta^{-1}$, then the
bistellar move $\beta$ is said to be {\it inessential}. Let a
bistellar move be {\it essential} if it is not inessential.

Let us construct the graph $\Gamma_n$ in the following way. The
vertex set of $\Gamma_n$ is the set $\mathcal T_{n+1}$. To each
equivalence class $\beta$ of essential bistellar moves transforming $L_1$
to $L_2$ assign the edge $e_{\beta}$ of $\Gamma_n$ with endpoints
$L_1$ and $L_2$. The edge $e_{\beta^{-1}}$ corresponding to the
inverse bistellar move coincides with the edge $e_{\beta}$ but
has the opposite orientation. By Pachner's theorem the graph
$\Gamma_n$ is connected.

Let $C_*(\Gamma_n;\mathbb{Z})$ be the cellular chain complex of
the graph $\Gamma_n$. We mean that $e_{\beta}=0\in
C_1(\Gamma_n;\mathbb{Z})$ if the bistellar move $\beta$ is
inessential. The group $\mathbb Z_2$ acts on the graph $\Gamma_n$
by changing orientation of all PL spheres. The group $\mathbb
Z_2$ acts on the group $\mathbb{Q}$ by changing sign. Let
$C^*_{\mathbb Z_2}(\Gamma_n;\mathbb{Q})=
\Hom_{\mathbb{Z}_2}(C_*(\Gamma_n;\mathbb{Z}),\mathbb{Q})$ be the
equivariant cochain complex of $\Gamma_n$. Let
$H^*_{\mathbb{Z}_2}(\Gamma_n;\mathbb{Q})=
H^*(C^*_{\mathbb{Z}_2}(\Gamma_n;\mathbb{Q}))$ be the equivariant
cohomology of $\Gamma_n$. By $d$ denote the differential of the
complex $C^*_{\mathbb Z_2}(\Gamma_n;\mathbb Q)$. We have the
isomorphism $H^1_{\mathbb Z_2}(\Gamma_n;\mathbb
Q)\cong\Hom_{\mathbb Z_2} (H_1(\Gamma_n;\mathbb Z),\mathbb Q)$.

Evidently, $C^0_{\mathbb{Z}_2}(\Gamma_{n-1};\mathbb{Q})=\mathcal
T^n(\mathbb Q)$. Hence we have the differential $\delta :
C^0_{\mathbb{Z}_2}(\Gamma_{n-1};\mathbb{Q})\to
C^0_{\mathbb{Z}_2}(\Gamma_n;\mathbb{Q})$, $\delta^2=0$. Consider
a bistellar move $\beta$ transforming $L_1$ to $L_2$, 
$L_1,L_2\in\mathcal T_{n+1}$. We may assume that $L_1$ and $L_2$ are
simplicial complexes on the same vertex set $V$. It is actually
true if $\beta$ is associated with the simplex, whose dimension is
neither $0$ nor $n$. Otherwise we assume that one of the two
simplicial complexes contains the vertex $v_0$ which is not a
simplex of this complex. Let $W$ be the set of all vertices $v\in
V$ such that the bistellar move $\beta$ induces an essential
bistellar move $\beta_v$ transforming $\Lk_{K_1}v$ to $\Lk_{K_2}v$
(by definition, $v_0\notin W$). Let the differential $\delta :
C^1_{\mathbb{Z}_2}(\Gamma_{n-1};\mathbb{Q})\to
C^1_{\mathbb{Z}_2}(\Gamma_n;\mathbb{Q})$ be given by $(\delta
h)(e_{\beta})=\sum_{v\in W} h(e_{\beta_v})$. It is easy to prove
that $\delta^2=0$ and $\delta d=d\delta$.

Put $C^{j,n}=C^j_{\mathbb{Z}_2}(\Gamma_{n-1};\mathbb{Q})$. Then
$C^{*,*}$ is a bigraded complex. We have $\bideg d=(1,0)$ and
$\bideg\delta=(0,1)$. By $Z^{*,*}_d$, $B^{*,*}_d$, and
$H^{*,*}_d$ denote respectively the cocycle group, the coboundary
group, and the cohomology group of the complex  $C^{*,*}$ with
respect to the differential $d$. By $Z^{*,*}_{\delta}$,
$B^{*,*}_{\delta}$ и $H^{*,*}_{\delta}$  denote respectively the
cocycle group, the coboundary group, and the cohomology group of
the complex  $C^{*,*}$ with respect to the differential $\delta$.
The graph $\Gamma_{n-1}$ is connected. Hence $H_d^{0,n}=0$.
Therefore $d:C^{0,n}\to C^{1,n}$ is a monomorphism.

Consider $L_1,L_2\in\mathcal T_{n}$. Let $\beta$, $V$, and $W$ be
as above. Consider the cone $CL_1$ with vertex $u_1$ and the cone
$CL_2$ with vertex $u_2$. Let $L_{\beta}$ be the simplicial
complex on the set $V\cup\{ u_1,u_2\}$ given by
$L_{\beta}=CL_1\cup CL_2\cup (\Delta_1 *\Delta_2)$. Choose the
orientation of $L_{\beta}$ such that the induced orientation of
$\Lk u_2$ coincides with the orientation of $L_2$. Then
$L_{\beta}\in\mathcal{T}_{n+1}$.  If $\beta_1$ and $\beta_2$
are equivalent bistellar moves, then the complexes $L_{\beta_1}$
and $L_{\beta_2}$ are isomorphic. The complexes $L_{\beta}$ and
$L_{\beta^{-1}}$ are antiisomorphic. If $\beta$ is inessential, then
$L_{\beta}$ is symmetric. Let the homomorphism $s:C^{0,n}\to
C^{1,n-1}$ be given by $s(f)(e_{\beta}) = f(L_{\beta})$.

Obviously, $d: C^{0,*}\to C^{1,*}$ is a chain homomorphism of the
complexes with differential $\delta$.

\begin{propos}
The homomorphism $s$ is a chain homotopy between the chain
homomorphisms  $d$ and $0$ of $C^{0,*}$ to $C^{1,*}$, i. e.
$d=\delta s +s\delta$.
\end{propos}

\begin{proof}
For any $v\in V\backslash W$ the link of the vertex $v$ in
$L_{\beta}$ is symmetric. For any $v\in W$ the link of the vertex
$v$ in $L_{\beta}$ is isomorphic to $-L_{\beta_v}$. The links of
the vertices $u_1$ and $u_2$ are isomorphic to $-L_1$ and $L_2$
respectively. Hence for any $f\in C^{0,n}$ we have

\begin{multline*}
s(\delta f)(e_{\beta})=(\delta f)(L_{\beta})=-\sum_{v\in
W}f(L_{\beta_v})+f(L_2)-f(L_1)\\= -\sum_{v\in
W}s(f)(e_{\beta_v})+f(\partial e_{\beta})=-\delta
s(f)(e_{\beta})+df(e_{\beta})
\end{multline*}

Consequently $df=\delta s(f)+s(\delta f)$.
\end{proof}

Let $A^n$ be the subgroup of $C^{1,n}$ consisting of all $h$ such
that $\delta h\in B^{1,n+1}_d$.

\begin{propos}
The homomorphism $s|_{Z^{0,n}_{\delta}}$ is a monomorphism and
$s(Z^{0,n}_{\delta})\subset A^{n-1}$.
\end{propos}

\begin{proof}
By Proposition 3.1, $d|_{Z^{0,n}_{\delta}}=\delta
s|_{Z^{0,n}_{\delta}}$. The homomorphism $d$ is a monomorphism.
Therefore $s|_{Z^{0,n}_{\delta}}$ is a monomorphism. If $f\in
Z^{0,n}_{\delta}$, then $\delta s(f)=df$. Hence $s(f)\in A^{n-1}$.
\end{proof}

\begin{corr}
$Z^{0,3}_{\delta}=0$. Hence $H^{0,3}_{\delta}=0$.
\end{corr}

\begin{proof}
The graph $\Gamma_1$ is isomorphic to the graph with the vertex
set $\{ 3,4,5,6,\ldots\}$ such that for any $k$ there exists a
unique edge with the endpoints $k$ and $(k+1)$. The action of the
group $\mathbb{Z}_2$ is trivial. Therefore $C^{1,2}=0$. By
Proposition 3.2, there exists a monomorphism of $Z^{0,3}_{\delta}$
to $C^{1,2}$. Hence $Z^{0,3}_{\delta}=0$.
\end{proof}

\begin{propos} $s|_{B^{0,4}_{\delta}}$ is an isomorphism of $B^{0,4}_{\delta}$
to $B^{1,3}_d$.
\end{propos}

\begin{proof} We have $C^{1,2}=0$. Therefore $dg=s(\delta g)$
for any $g\in C^{0,3}$. The proposition follows.
\end{proof}

Thus $s$ induces the monomorphism $s^*$ of $H^{0,4}_{\delta}$ to
$H^{1,3}_d=H^1_{\mathbb Z_2}(\Gamma_2;\mathbb Q)$. Let $\widetilde
A^3$ be the kernel of the homomorphism $\delta^*:H^{1,3}\to
H^{1,4}$ induced by the chain homomorphism $\delta : C^{*,3}\to
C^{*,4}$. Then $s^*(H^{0,4}_{\delta})\subset \widetilde A^{3}$.

\subsection{Generators of the group $H_1(\Gamma_2;\mathbb{Z})$}

In the sequel if we say that $\{ u_0,\ldots,u_l\}$ is a simplex of
an $l$-dimensional oriented simplicial complex $L$, then we mean
that the sequence of the vertices $u_0,\ldots,u_l$ provides the
given orientation of $L$. If we show a $2$-dimensional simplicial
complex in a figure, then we mean that the orientation is
clockwise. For any cycle $\gamma\in Z_1(\Gamma_n;\mathbb{Z})$ let
$\bar\gamma\in H_1(\Gamma_n;\mathbb Z)$ be the homology class
represented by $\gamma$.

Let $L$ be an oriented 2-dimensional PL sphere. An edge $e\in L$
is called {\it admissible} if there exists a bistellar move
associated with $e$. A pair of edges $(e_1,e_2)$ is called {\it
admissible} if:\\
1) there is no triangle $\Delta\in K$ containing both $e_1$ and
$e_2$;\\
2) $e_1$ and $e_2$ are admissible;\\
3) the bistellar move associated with the edge $e_1$ takes $e_2$
to an admissible edge.

Let $\Delta_1,\Delta_2\in L$ be two distinct triangles. Apply to
$L$ the bistellar move associated with $\Delta_1$; by $v_1$ denote
the new vertex created by this bistellar move. Then apply to the
obtained triangulation the bistellar move associated with
$\Delta_2$; by $v_2$ denote the new vertex created by this
bistellar move. Apply  to the obtained triangulation the
bistellar move associated with the vertex $v_1$. Finally, apply
to the obtained triangulation  the bistellar move associated with
the vertex $v_2$. By $\alpha_1(L,\Delta_1,\Delta_2)$ denote the
obtained cycle in the graph $\Gamma_2$, see Fig.~1,~a,~b,~c.

\begin{picture}(230,193)
\put(0,137){
\begin{picture}(102,56)

\put(0,39){
     \begin{picture}(47,16)
     \put(1,1){\circle*{1}}
     \put(11,16){\circle*{1}}
     \put(21,1){\circle*{1}}
     \put(26,1){\circle*{1}}
     \put(36,16){\circle*{1}}
     \put(46,1){\circle*{1}}
     \put(1,1){\line(1,0){20}}
     \put(1,1){\line(2,3){10}}
     \put(11,16){\line(2,-3){10}}
     \put(26,1){\line(1,0){20}}
     \put(26,1){\line(2,3){10}}
     \put(36,16){\line(2,-3){10}}
     \end{picture}
    }
\put(55,39){
     \begin{picture}(47,16)
     \put(1,1){\circle*{1}}
     \put(11,16){\circle*{1}}
     \put(21,1){\circle*{1}}
     \put(26,1){\circle*{1}}
     \put(36,16){\circle*{1}}
     \put(46,1){\circle*{1}}
     \put(1,1){\line(1,0){20}}
     \put(1,1){\line(2,3){10}}
     \put(11,16){\line(2,-3){10}}
     \put(26,1){\line(1,0){20}}
     \put(26,1){\line(2,3){10}}
     \put(36,16){\line(2,-3){10}}
     \end{picture}
    }
\put(55,9){
     \begin{picture}(47,16)
     \put(1,1){\circle*{1}}
     \put(11,16){\circle*{1}}
     \put(21,1){\circle*{1}}
     \put(26,1){\circle*{1}}
     \put(36,16){\circle*{1}}
     \put(46,1){\circle*{1}}
     \put(1,1){\line(1,0){20}}
     \put(1,1){\line(2,3){10}}
     \put(11,16){\line(2,-3){10}}
     \put(26,1){\line(1,0){20}}
     \put(26,1){\line(2,3){10}}
     \put(36,16){\line(2,-3){10}}
     \end{picture}
    }
\put(0,9){
     \begin{picture}(47,16)
     \put(1,1){\circle*{1}}
     \put(11,16){\circle*{1}}
     \put(21,1){\circle*{1}}
     \put(26,1){\circle*{1}}
     \put(36,16){\circle*{1}}
     \put(46,1){\circle*{1}}
     \put(1,1){\line(1,0){20}}
     \put(1,1){\line(2,3){10}}
     \put(11,16){\line(2,-3){10}}
     \put(26,1){\line(1,0){20}}
     \put(26,1){\line(2,3){10}}
     \put(36,16){\line(2,-3){10}}
     \end{picture}
    }
\put(47.5,47){\vector(1,0){10}}
\put(81,35){\vector(0,-1){10}}
\put(57.5,17){\vector(-1,0){10}}
\put(26,25){\vector(0,1){10}}
\put(10,44){$\Delta_1$}
\put(35,44){$\Delta_2$}
\put(56,40){
    \begin{picture}(20,15)
    \put(10,5){\circle*{1}}
    \put(0,0){\line(2,1){10}}
    \put(10,5){\line(0,1){10}}
    \put(10,5){\line(2,-1){10}}
    \end{picture}
    }
\put(56,10){
    \begin{picture}(20,15)
    \put(10,5){\circle*{1}}
    \put(0,0){\line(2,1){10}}
    \put(10,5){\line(0,1){10}}
    \put(10,5){\line(2,-1){10}}
    \end{picture}
    }
\put(81,10){
    \begin{picture}(20,15)
    \put(10,5){\circle*{1}}
    \put(0,0){\line(2,1){10}}
    \put(10,5){\line(0,1){10}}
    \put(10,5){\line(2,-1){10}}
    \end{picture}
    }
\put(26,10){
    \begin{picture}(20,15)
    \put(10,5){\circle*{1}}
    \put(0,0){\line(2,1){10}}
    \put(10,5){\line(0,1){10}}
    \put(10,5){\line(2,-1){10}}
    \end{picture}
    }
\put(51,0){a}
\end{picture}

}
\put(0,76){
\begin{picture}(102,56)
\put(0,39){
     \begin{picture}(47,16)
     \put(1,1){\circle*{1}}
     \put(1,16){\circle*{1}}
     \put(19.75,8.5){\circle*{1}}
     \put(38.5,16){\circle*{1}}
     \put(38.5,1){\circle*{1}}
     \put(1,1){\line(0,1){15}}
     \put(1,1){\line(5,2){18.75}}
     \put(1,16){\line(5,-2){18.75}}
     \put(38.5,1){\line(0,1){15}}
     \put(38.5,1){\line(-5,2){18.75}}
     \put(38.5,16){\line(-5,-2){18.75}}
     \qbezier(23.25,9.9)(19.75,12)(16.25,9.9)
     \put(23.7,5){${}^x$}
     \put(18.3,12.5){$\vartheta_1$}
     \qbezier(23.25,7.1)(19.75,5)(16.25,7.1)
     \put(18.3,0){$\vartheta_2$}
     \end{picture}
    }
\put(62.5,39){
     \begin{picture}(47,16)
     \put(1,1){\circle*{1}}
     \put(1,16){\circle*{1}}
     \put(19.75,8.5){\circle*{1}}
     \put(38.5,16){\circle*{1}}
     \put(38.5,1){\circle*{1}}
     \put(4.75,8.5){\circle*{1}}
     \put(1,1){\line(0,1){15}}
     \put(1,1){\line(5,2){18.75}}
     \put(1,16){\line(5,-2){18.75}}
     \put(38.5,1){\line(0,1){15}}
     \put(38.5,1){\line(-5,2){18.75}}
     \put(38.5,16){\line(-5,-2){18.75}}
     \put(1,1){\line(1,2){3.75}}
     \put(1,16){\line(1,-2){3.75}}
     \put(4.75,8.5){\line(1,0){15}}
     \end{picture}
    }
\put(62.5,9){
     \begin{picture}(47,16)
     \put(1,1){\circle*{1}}
     \put(1,16){\circle*{1}}
     \put(19.75,8.5){\circle*{1}}
     \put(38.5,16){\circle*{1}}
     \put(38.5,1){\circle*{1}}
     \put(4.75,8.5){\circle*{1}}
     \put(34.75,8.5){\circle*{1}}
     \put(1,1){\line(0,1){15}}
     \put(1,1){\line(5,2){18.75}}
     \put(1,16){\line(5,-2){18.75}}
     \put(38.5,1){\line(0,1){15}}
     \put(38.5,1){\line(-5,2){18.75}}
     \put(38.5,16){\line(-5,-2){18.75}}
     \put(1,1){\line(1,2){3.75}}
     \put(1,16){\line(1,-2){3.75}}
     \put(4.75,8.5){\line(1,0){15}}
     \put(38.5,1){\line(-1,2){3.75}}
     \put(38.5,16){\line(-1,-2){3.75}}
     \put(19.75,8.5){\line(1,0){15}}
     \end{picture}
    }
\put(0,9){
     \begin{picture}(47,16)
     \put(1,1){\circle*{1}}
     \put(1,16){\circle*{1}}
     \put(19.75,8.5){\circle*{1}}
     \put(38.5,16){\circle*{1}}
     \put(38.5,1){\circle*{1}}
     \put(34.75,8.5){\circle*{1}}
     \put(1,1){\line(0,1){15}}
     \put(1,1){\line(5,2){18.75}}
     \put(1,16){\line(5,-2){18.75}}
     \put(38.5,1){\line(0,1){15}}
     \put(38.5,1){\line(-5,2){18.75}}
     \put(38.5,16){\line(-5,-2){18.75}}
     \put(38.5,1){\line(-1,2){3.75}}
     \put(38.5,16){\line(-1,-2){3.75}}
     \put(19.75,8.5){\line(1,0){15}}
     \end{picture}
    }
\put(47.5,47){\vector(1,0){10}}
\put(84.5,35){\vector(0,-1){10}}
\put(57.5,17){\vector(-1,0){10}}
\put(22,25){\vector(0,1){10}}
\put(5,46){$\Delta_1$}
\put(31,46){$\Delta_2$}

\put(51,0){b}
\end{picture}

}
\put(0,10){
\begin{picture}(102,61)

\put(0,41){
     \begin{picture}(47,16)
     \put(1,8.5){\circle*{1}}
     \put(19.75,16){\circle*{1}}
     \put(19.75,1){\circle*{1}}
     \put(38.5,8.5){\circle*{1}}
     \put(19.75,1){\line(0,1){15}}
     \put(1,8.5){\line(5,2){18.75}}
     \put(1,8.5){\line(5,-2){18.75}}
     \put(38.5,8.5){\line(-5,2){18.75}}
     \put(38.5,8.5){\line(-5,-2){18.75}}
     \put(18,17){$x$}
     \put(18,-3){$y$}
     \end{picture}
    } 
\put(62.5,41){
     \begin{picture}(47,16)
     \put(1,8.5){\circle*{1}}
     \put(16,8.5){\circle*{1}}
     \put(19.75,16){\circle*{1}}
     \put(19.75,1){\circle*{1}}
     \put(38.5,8.5){\circle*{1}}
     \put(19.75,1){\line(0,1){15}}
     \put(1,8.5){\line(5,2){18.75}}
     \put(1,8.5){\line(5,-2){18.75}}
     \put(38.5,8.5){\line(-5,2){18.75}}
     \put(38.5,8.5){\line(-5,-2){18.75}}
     \put(16,8.5){\line(1,2){3.75}}
     \put(16,8.5){\line(1,-2){3.75}}
     \put(1,8.5){\line(1,0){15}}
     \end{picture}
    } 
\put(62.5,9){
     \begin{picture}(47,16)
     \put(1,8.5){\circle*{1}}
     \put(16,8.5){\circle*{1}}
     \put(23.5,8.5){\circle*{1}}
     \put(19.75,16){\circle*{1}}
     \put(19.75,1){\circle*{1}}
     \put(38.5,8.5){\circle*{1}}
     \put(19.75,1){\line(0,1){15}}
     \put(1,8.5){\line(5,2){18.75}}
     \put(1,8.5){\line(5,-2){18.75}}
     \put(38.5,8.5){\line(-5,2){18.75}}
     \put(38.5,8.5){\line(-5,-2){18.75}}
     \put(16,8.5){\line(1,2){3.75}}
     \put(16,8.5){\line(1,-2){3.75}}
     \put(1,8.5){\line(1,0){15}}
     \put(23.5,8.5){\line(-1,2){3.75}}
     \put(23.5,8.5){\line(-1,-2){3.75}}
     \put(23.5,8.5){\line(1,0){15}}
     \end{picture}
    } 
\put(0,9){
     \begin{picture}(47,16)
     \put(1,8.5){\circle*{1}}
     \put(23.5,8.5){\circle*{1}}
     \put(19.75,16){\circle*{1}}
     \put(19.75,1){\circle*{1}}
     \put(38.5,8.5){\circle*{1}}
     \put(19.75,1){\line(0,1){15}}
     \put(1,8.5){\line(5,2){18.75}}
     \put(1,8.5){\line(5,-2){18.75}}
     \put(38.5,8.5){\line(-5,2){18.75}}
     \put(38.5,8.5){\line(-5,-2){18.75}}
     
     \put(23.5,8.5){\line(-1,2){3.75}}
     \put(23.5,8.5){\line(-1,-2){3.75}}
     \put(23.5,8.5){\line(1,0){15}}
     \end{picture}
    } 
\put(47.5,49){\vector(1,0){10}}
\put(84.5,38.5){\vector(0,-1){10}}
\put(57.5,17){\vector(-1,0){10}}
\put(22,26){\vector(0,1){10}}
\put(12,48){$\Delta_1$}
\put(25,48){$\Delta_2$}

\put(51,0){c}
\end{picture}

}
\put(128,137){
\begin{picture}(102,56)

\put(0,39){
     \begin{picture}(47,16)
     \put(1,1){\circle*{1}}
     \put(11,16){\circle*{1}}
     \put(21,1){\circle*{1}}
     \put(26,1){\circle*{1}}
     \put(26,16){\circle*{1}}
     \put(41,1){\circle*{1}}
     \put(41,16){\circle*{1}}
     \put(1,1){\line(1,0){20}}
     \put(1,1){\line(2,3){10}}
     \put(11,16){\line(2,-3){10}}
     \put(26,1){\line(1,0){15}}
     \put(26,1){\line(0,1){15}}
     \put(26,16){\line(1,0){15}}
     \put(41,1){\line(0,1){15}}
     \put(26,16){\line(1,-1){15}}
     \end{picture}
    } 
\put(55,39){
     \begin{picture}(47,16)
     \put(1,1){\circle*{1}}
     \put(11,16){\circle*{1}}
     \put(21,1){\circle*{1}}
     \put(26,1){\circle*{1}}
     \put(26,16){\circle*{1}}
     \put(41,1){\circle*{1}}
     \put(41,16){\circle*{1}}
     \put(1,1){\line(1,0){20}}
     \put(1,1){\line(2,3){10}}
     \put(11,16){\line(2,-3){10}}
     \put(26,1){\line(1,0){15}}
     \put(26,1){\line(0,1){15}}
     \put(26,16){\line(1,0){15}}
     \put(41,1){\line(0,1){15}}
     \put(26,16){\line(1,-1){15}}
     \end{picture}
    } 
\put(55,9){
     \begin{picture}(47,16)
     \put(1,1){\circle*{1}}
     \put(11,16){\circle*{1}}
     \put(21,1){\circle*{1}}
     \put(26,1){\circle*{1}}
     \put(26,16){\circle*{1}}
     \put(41,1){\circle*{1}}
     \put(41,16){\circle*{1}}
     \put(1,1){\line(1,0){20}}
     \put(1,1){\line(2,3){10}}
     \put(11,16){\line(2,-3){10}}
     \put(26,1){\line(1,0){15}}
     \put(26,1){\line(0,1){15}}
     \put(26,16){\line(1,0){15}}
     \put(41,1){\line(0,1){15}}
     \put(26,1){\line(1,1){15}}
     \end{picture}
    } 
\put(0,9){
     \begin{picture}(47,16)
     \put(1,1){\circle*{1}}
     \put(11,16){\circle*{1}}
     \put(21,1){\circle*{1}}
     \put(26,1){\circle*{1}}
     \put(26,16){\circle*{1}}
     \put(41,1){\circle*{1}}
     \put(41,16){\circle*{1}}
     \put(1,1){\line(1,0){20}}
     \put(1,1){\line(2,3){10}}
     \put(11,16){\line(2,-3){10}}
     \put(26,1){\line(1,0){15}}
     \put(26,1){\line(0,1){15}}
     \put(26,16){\line(1,0){15}}
     \put(41,1){\line(0,1){15}}
     \put(26,1){\line(1,1){15}}
     \end{picture}
    } 
\put(46,47){\vector(1,0){10}}
\put(77,35){\vector(0,-1){10}}
\put(56,17){\vector(-1,0){10}}
\put(23,25){\vector(0,1){10}}
\put(10,44){$\Delta$}
\put(32,45){$e$}
\put(56,40){
    \begin{picture}(20,15)
    \put(10,5){\circle*{1}}
    \put(0,0){\line(2,1){10}}
    \put(10,5){\line(0,1){10}}
    \put(10,5){\line(2,-1){10}}
    \end{picture}
    }
\put(56,10){
    \begin{picture}(20,15)
    \put(10,5){\circle*{1}}
    \put(0,0){\line(2,1){10}}
    \put(10,5){\line(0,1){10}}
    \put(10,5){\line(2,-1){10}}
    \end{picture}
    }

\put(51,0){a}
\end{picture}

}
\put(128,76){
\begin{picture}(102,60)

\put(0,39){
     \begin{picture}(47,16)
     \put(1,1){\circle*{1}}
     \put(1,16){\circle*{1}}
     \put(19.75,8.5){\circle*{1}}
     \put(43.75,8.5){\circle*{1}}
     \put(31.75,20.5){\circle*{1}}
     \put(31.75,-3.5){\circle*{1}}
     \put(1,1){\line(0,1){15}}
     \put(1,1){\line(5,2){18.75}}
     \put(1,16){\line(5,-2){18.75}}
     \put(31.75,-3.5){\line(0,1){16}}
     \put(31.75,20.5){\line(0,-1){4}}
     \put(19.75,8.5){\line(1,1){12}}
     \put(19.75,8.5){\line(1,-1){12}}
     \put(31.75,-3.5){\line(1,1){12}}
     \put(31.75,20.5){\line(1,-1){12}}
     \qbezier(22.25,11.25)(17.5,12.2)(16.25,9.9)
     \put(17,13.5){$\vartheta_1$}
     \qbezier(22.25,5.75)(17.5,4.7)(16.25,7.1)
     \put(17,0){$\vartheta_2$}
     \end{picture}
    } 
\put(55,39){
     \begin{picture}(47,16)
     \put(1,1){\circle*{1}}
     \put(1,16){\circle*{1}}
     \put(19.75,8.5){\circle*{1}}
     \put(43.75,8.5){\circle*{1}}
     \put(31.75,20.5){\circle*{1}}
     \put(31.75,-3.5){\circle*{1}}
     \put(4.75,8.5){\circle*{1}}
     \put(1,1){\line(0,1){15}}
     \put(1,1){\line(5,2){18.75}}
     \put(1,16){\line(5,-2){18.75}}
     \put(31.75,-3.5){\line(0,1){24}}
     \put(19.75,8.5){\line(1,1){12}}
     \put(19.75,8.5){\line(1,-1){12}}
     \put(31.75,-3.5){\line(1,1){12}}
     \put(31.75,20.5){\line(1,-1){12}}
     \put(4.75,8.5){\line(1,0){15}}
     \put(1,1){\line(1,2){3.75}}
     \put(1,16){\line(1,-2){3.75}}
     \end{picture}
    } 
\put(55,9){
     \begin{picture}(47,16)
     \put(1,1){\circle*{1}}
     \put(1,16){\circle*{1}}
     \put(19.75,8.5){\circle*{1}}
     \put(43.75,8.5){\circle*{1}}
     \put(31.75,20.5){\circle*{1}}
     \put(31.75,-3.5){\circle*{1}}
     \put(4.75,8.5){\circle*{1}}
     \put(1,1){\line(0,1){15}}
     \put(1,1){\line(5,2){18.75}}
     \put(1,16){\line(5,-2){18.75}}
     \put(19.75,8.5){\line(1,0){24}}
     \put(19.75,8.5){\line(1,1){12}}
     \put(19.75,8.5){\line(1,-1){12}}
     \put(31.75,-3.5){\line(1,1){12}}
     \put(31.75,20.5){\line(1,-1){12}}
     \put(4.75,8.5){\line(1,0){15}}
     \put(1,1){\line(1,2){3.75}}
     \put(1,16){\line(1,-2){3.75}}
     \end{picture}
    } 
\put(0,9){
     \begin{picture}(47,16)
     \put(1,1){\circle*{1}}
     \put(1,16){\circle*{1}}
     \put(19.75,8.5){\circle*{1}}
     \put(43.75,8.5){\circle*{1}}
     \put(31.75,20.5){\circle*{1}}
     \put(31.75,-3.5){\circle*{1}}
     
     \put(1,1){\line(0,1){15}}
     \put(1,1){\line(5,2){18.75}}
     \put(1,16){\line(5,-2){18.75}}
     \put(19.75,8.5){\line(1,0){24}}
     \put(19.75,8.5){\line(1,1){12}}
     \put(19.75,8.5){\line(1,-1){12}}
     \put(31.75,-3.5){\line(1,1){12}}
     \put(31.75,20.5){\line(1,-1){12}}
     \end{picture}
    } 
\put(47.3,47){\vector(1,0){9}}
\put(77,35){\vector(0,-1){10}}
\put(56.3,17){\vector(-1,0){9}}
\put(23,25){\vector(0,1){10}}
\put(6,45){$\Delta$}
\put(33,52){$e$}
\put(25,45){$\Delta_1$}
\put(35,45){$\Delta_2$}

\put(51,0){b}
\end{picture}

}
\put(128,10){
\begin{picture}(102,61)

\put(0,41){
     \begin{picture}(47,16)
     \put(1,8.5){\circle*{1}}
     \put(19.75,16){\circle*{1}}
     \put(19.75,1){\circle*{1}}
     \put(34.75,1){\circle*{1}}
     \put(34.75,16){\circle*{1}}
     \put(19.75,1){\line(0,1){15}}
     \put(1,8.5){\line(5,2){18.75}}
     \put(1,8.5){\line(5,-2){18.75}}
     \put(19.75,1){\line(1,0){15}}
     \put(19.75,16){\line(1,0){15}}
     \put(34.75,1){\line(0,1){15}}
     \put(19.75,16){\line(1,-1){15}}
     
     \end{picture}
    } 
\put(61.25,41){
     \begin{picture}(47,16)
     \put(1,8.5){\circle*{1}}
     \put(19.75,16){\circle*{1}}
     \put(19.75,1){\circle*{1}}
     \put(34.75,1){\circle*{1}}
     \put(34.75,16){\circle*{1}}
     \put(16,8.5){\circle*{1}}
     \put(19.75,1){\line(0,1){15}}
     \put(1,8.5){\line(5,2){18.75}}
     \put(1,8.5){\line(5,-2){18.75}}
     \put(19.75,1){\line(1,0){15}}
     \put(19.75,16){\line(1,0){15}}
     \put(34.75,1){\line(0,1){15}}
     \put(19.75,16){\line(1,-1){15}}
     \put(16,8.5){\line(1,2){3.75}}
     \put(16,8.5){\line(1,-2){3.75}}
     \put(1,8.5){\line(1,0){15}}
     \end{picture}
    } 
\put(61.25,9){
     \begin{picture}(47,16)
     \put(1,8.5){\circle*{1}}
     \put(19.75,16){\circle*{1}}
     \put(19.75,1){\circle*{1}}
     \put(34.75,1){\circle*{1}}
     \put(34.75,16){\circle*{1}}
     \put(16,8.5){\circle*{1}}
     \put(19.75,1){\line(0,1){15}}
     \put(1,8.5){\line(5,2){18.75}}
     \put(1,8.5){\line(5,-2){18.75}}
     \put(19.75,1){\line(1,0){15}}
     \put(19.75,16){\line(1,0){15}}
     \put(34.75,1){\line(0,1){15}}
     \put(19.75,1){\line(1,1){15}}
     \put(16,8.5){\line(1,2){3.75}}
     \put(16,8.5){\line(1,-2){3.75}}
     \put(1,8.5){\line(1,0){15}}
     \end{picture}
    } 
\put(0,9){
     \begin{picture}(47,16)
     \put(1,8.5){\circle*{1}}
     \put(19.75,16){\circle*{1}}
     \put(19.75,1){\circle*{1}}
     \put(34.75,1){\circle*{1}}
     \put(34.75,16){\circle*{1}}
     \put(19.75,1){\line(0,1){15}}
     \put(1,8.5){\line(5,2){18.75}}
     \put(1,8.5){\line(5,-2){18.75}}
     \put(19.75,1){\line(1,0){15}}
     \put(19.75,16){\line(1,0){15}}
     \put(34.75,1){\line(0,1){15}}
     \put(19.75,1){\line(1,1){15}}
     \end{picture}
    } 
\put(46,47){\vector(1,0){10}}
\put(82,38){\vector(0,-1){10}}
\put(56,17){\vector(-1,0){10}}
\put(23,26){\vector(0,1){10}}
\put(13.5,47){$\Delta$}
\put(32.2,46.9){$e$}
\put(28,51.6){$\Delta_2$}
\put(22.5,45.2){$\Delta_1$}
\put(19,37){$x$}
\put(19,59){$y$}

\put(51,0){c}
\end{picture}

}
\put(45,0){Fig. 1}
\put(173,0){Fig. 2}
\end{picture}

There are three possibilities: the triangles $\Delta_1$ and
$\Delta_2$ can have no common vertices or have $1$ or $2$ common
vertices. Let $\mathcal S_1^0$ be the set of all homology classes
$\bar\alpha_1(L,\Delta_1,\Delta_2)\in H_1(\Gamma_2,\mathbb Z)$
such that the triangles $\Delta_1$ and $\Delta_2$  have no common
vertices, see Fig. 1, a.

By $\mathcal S_1^1(p,q)$ denote the set of all homology classes
$\bar\alpha_1(L,\Delta_1,\Delta_2)$ such that:\\
1) the triangles $\Delta_1$ and $\Delta_2$ have a unique common
vertex $x$, see Fig. 1, b;\\
2) there are exactly $p$ triangles containing the vertex $x$ and
situated in the angle~$\vartheta_1$;\\
3) there are exactly $q$ triangles containing the vertex $x$ and
situated in the angle~$\vartheta_2$.

By $\mathcal S_1^2(p,q)$ denote the set of all homology classes
$\bar\alpha_1(L,\Delta_1,\Delta_2)$ such that:\\
1) the triangles $\Delta_1$ and $\Delta_2$ have a  common
edge $e$ with endpoints $x$ and $y$, see Fig. 1, c;\\
2) there are exactly $p$ triangles that contain the vertex $x$
and coincide neither with $\Delta_1$ nor with $\Delta_2$;\\
3) there are exactly $q$ triangles that contain the vertex $y$ and
coincide neither with $\Delta_1$ nor with $\Delta_2$.

Let $\Delta\in L$ be a triangle. Let $e\in L$ be an admissible
edge such that $e\not\subset\Delta$. By $\alpha_2(L,\Delta,e)$
denote the cycle shown in Fig.~2,~a,~b,~c. Let $\mathcal S_2^0$ be
the set of all homology classes $\bar\alpha_2(L,\Delta,e)$ such
that the triangle $\Delta$ has common vertices neither with
$\Delta_1$ nor with $\Delta_2$, where $\Delta_1$ and $\Delta_2$
are the two triangles containing the edge $e$, see Fig.~2,~a.

By $\mathcal S_2^1(p,q)$ denote the set of all homology classes
$\bar\alpha_2(L,\Delta,e)$ such that:\\
1) the triangle $\Delta$ has a unique common vertex $x$ with the
triangle $\Delta_1$, see Fig.~2,~b;\\
2) there are exactly $p$ triangles containing the vertex $x$ and
situated in the angle~$\vartheta_1$;\\
3) there are exactly $q$ triangles containing the vertex $x$ and
situated in the angle~$\vartheta_2$.

By $\mathcal S_2^2(p,q)$ denote the set of all homology classes
$\bar\alpha_2(L,\Delta,e)$ such that:\\
1) the triangles $\Delta$ and $\Delta_1$ have a  common
edge $e_1$ with endpoints $x$ and $y$, see Fig. 2, c;\\
2) there are exactly $p$ triangles that contain the vertex $x$
and coincide neither with $\Delta$ nor with $\Delta_1$;\\
3) there are exactly $q$ triangles that contain the vertex $y$ and
coincide with none of the triangles $\Delta$, $\Delta_1$, and
$\Delta_2$.

Assume that $(e_1,e_2)$ is an admissible pair. By
$\alpha_3(L,e_1,e_2)$ denote the cycle shown in Fig.~3,~a,~b,~c.
Let $\mathcal S_3^0$ be the set of all homology classes
$\bar\alpha_3(L,e_1,e_2)$ such that any triangle containing the
edge $e_1$ has no common vertices  with any triangle containing
the edge $e_2$, see Fig.~3,~a.

\begin{picture}(230,143)
\put(0,85){

\begin{picture}(97,56)
\put(-5,39){
     \begin{picture}(47,16)
     \put(6,1){\circle*{1}}
     \put(6,16){\circle*{1}}
     \put(21,1){\circle*{1}}
     \put(21,16){\circle*{1}}
     \put(26,1){\circle*{1}}
     \put(26,16){\circle*{1}}
     \put(41,1){\circle*{1}}
     \put(41,16){\circle*{1}}
     \put(6,1){\line(1,0){15}}
     \put(6,16){\line(1,0){15}}
     \put(6,1){\line(0,1){15}}
     \put(21,1){\line(0,1){15}}
     \put(6,16){\line(1,-1){15}}
     \put(26,1){\line(1,0){15}}
     \put(26,1){\line(0,1){15}}
     \put(26,16){\line(1,0){15}}
     \put(41,1){\line(0,1){15}}
     \put(26,16){\line(1,-1){15}}
     \end{picture}
    } 
\put(50,39){
     \begin{picture}(47,16)

     \put(6,1){\circle*{1}}
     \put(6,16){\circle*{1}}
     \put(21,1){\circle*{1}}
     \put(21,16){\circle*{1}}
     \put(26,1){\circle*{1}}
     \put(26,16){\circle*{1}}
     \put(41,1){\circle*{1}}
     \put(41,16){\circle*{1}}
     \put(6,1){\line(1,0){15}}
     \put(6,16){\line(1,0){15}}
     \put(6,1){\line(0,1){15}}
     \put(21,1){\line(0,1){15}}
     \put(6,1){\line(1,1){15}}
     \put(26,1){\line(1,0){15}}
     \put(26,1){\line(0,1){15}}
     \put(26,16){\line(1,0){15}}
     \put(41,1){\line(0,1){15}}
     \put(26,16){\line(1,-1){15}}     
     \end{picture}
    } 
\put(50,9){
     \begin{picture}(47,16)

     \put(6,1){\circle*{1}}
     \put(6,16){\circle*{1}}
     \put(21,1){\circle*{1}}
     \put(21,16){\circle*{1}}
     \put(26,1){\circle*{1}}
     \put(26,16){\circle*{1}}
     \put(41,1){\circle*{1}}
     \put(41,16){\circle*{1}}
     \put(6,1){\line(1,0){15}}
     \put(6,16){\line(1,0){15}}
     \put(6,1){\line(0,1){15}}
     \put(21,1){\line(0,1){15}}
     \put(6,1){\line(1,1){15}}
     \put(26,1){\line(1,0){15}}
     \put(26,1){\line(0,1){15}}
     \put(26,16){\line(1,0){15}}
     \put(41,1){\line(0,1){15}}
     \put(26,1){\line(1,1){15}}     
     \end{picture}
    } 
\put(-5,9){
     \begin{picture}(47,16)

     \put(6,1){\circle*{1}}
     \put(6,16){\circle*{1}}
     \put(21,1){\circle*{1}}
     \put(21,16){\circle*{1}}
     \put(26,1){\circle*{1}}
     \put(26,16){\circle*{1}}
     \put(41,1){\circle*{1}}
     \put(41,16){\circle*{1}}
     \put(6,1){\line(1,0){15}}
     \put(6,16){\line(1,0){15}}
     \put(6,1){\line(0,1){15}}
     \put(21,1){\line(0,1){15}}
     \put(6,16){\line(1,-1){15}}
     \put(26,1){\line(1,0){15}}
     \put(26,1){\line(0,1){15}}
     \put(26,16){\line(1,0){15}}
     \put(41,1){\line(0,1){15}}
     \put(26,1){\line(1,1){15}}     
     \end{picture}
    } 
\put(41,47){\vector(1,0){10}}
\put(75,38){\vector(0,-1){10}}
\put(51,17){\vector(-1,0){10}}
\put(20,28){\vector(0,1){10}}
\put(7,44){$e_1$}
\put(27,44){$e_2$}

\put(46,0){a}
\end{picture}

}
\put(115,79){
\begin{picture}(102,64)

\put(0,39){
     \begin{picture}(47,16)
     \put(1,13){\circle*{1}}
     \put(25,13){\circle*{1}}
     \put(13,1){\circle*{1}}
     \put(13,25){\circle*{1}}
     \put(37,1){\circle*{1}}
     \put(37,25){\circle*{1}}
     \put(49,13){\circle*{1}}
     \put(1,13){\line(1,1){12}}
     \put(1,13){\line(1,-1){12}}
     \put(13,1){\line(1,1){24}}
     \put(13,25){\line(1,-1){24}}
     \put(49,13){\line(-1,-1){12}}
     \put(49,13){\line(-1,1){12}} 
     \put(13,1){\line(0,1){16}}
     \put(37,1){\line(0,1){16}}
     \put(13,25){\line(0,-1){4}}
     \put(37,25){\line(0,-1){4}}
     \qbezier(27.5,15.5)(25,18)(22.5,15.5)
     \put(23,18){$\vartheta_1$}
     \qbezier(27.5,10.5)(25,8)(22.5,10.5)
     \put(23,3){$\vartheta_2$}
     \end{picture}
    } 
\put(62,39){
     \begin{picture}(47,16)

     \put(1,13){\circle*{1}}
     \put(25,13){\circle*{1}}
     \put(13,1){\circle*{1}}
     \put(13,25){\circle*{1}}
     \put(37,1){\circle*{1}}
     \put(37,25){\circle*{1}}
     \put(49,13){\circle*{1}}
     \put(1,13){\line(1,1){12}}
     \put(1,13){\line(1,-1){12}}
     \put(13,1){\line(1,1){24}}
     \put(13,25){\line(1,-1){24}}
     \put(49,13){\line(-1,-1){12}}
     \put(49,13){\line(-1,1){12}} 
     \put(1,13){\line(1,0){24}}
     \put(37,1){\line(0,1){24}}
     \end{picture}
    } 
\put(62,9){
     \begin{picture}(47,16)

     \put(1,13){\circle*{1}}
     \put(25,13){\circle*{1}}
     \put(13,1){\circle*{1}}
     \put(13,25){\circle*{1}}
     \put(37,1){\circle*{1}}
     \put(37,25){\circle*{1}}
     \put(49,13){\circle*{1}}
     \put(1,13){\line(1,1){12}}
     \put(1,13){\line(1,-1){12}}
     \put(13,1){\line(1,1){24}}
     \put(13,25){\line(1,-1){24}}
     \put(49,13){\line(-1,-1){12}}
     \put(49,13){\line(-1,1){12}} 
     \put(1,13){\line(1,0){48}}
    
     \end{picture}
    } 
\put(0,9){
     \begin{picture}(47,16)

     \put(1,13){\circle*{1}}
     \put(25,13){\circle*{1}}
     \put(13,1){\circle*{1}}
     \put(13,25){\circle*{1}}
     \put(37,1){\circle*{1}}
     \put(37,25){\circle*{1}}
     \put(49,13){\circle*{1}}
     \put(1,13){\line(1,1){12}}
     \put(1,13){\line(1,-1){12}}
     \put(13,1){\line(1,1){24}}
     \put(13,25){\line(1,-1){24}}
     \put(49,13){\line(-1,-1){12}}
     \put(49,13){\line(-1,1){12}} 
     \put(13,1){\line(0,1){24}}
     \put(25,13){\line(1,0){24}} 
     \end{picture}
    } 
\put(53.5,51.5){\vector(1,0){10}}
\put(89.5,42){\vector(0,-1){10}}
\put(63.5,21.5){\vector(-1,0){10}}
\put(27.5,28){\vector(0,1){10}}
\put(14,56){$e_1$}
\put(38,56){$e_2$}
\put(16,49){$\Delta_1$}
\put(30,49){$\Delta_2$}
\put(6,49){$\Delta_3$}
\put(40,49){$\Delta_4$}

\put(57,0){b}
\end{picture}

}
\put(50,10){
\begin{picture}(102,64)

\put(0,46){
     \begin{picture}(47,16)
     \put(1,1){\circle*{1}}
     \put(1,21){\circle*{1}}
     \put(21,1){\circle*{1}}
     \put(21,21){\circle*{1}}
     \put(41,1){\circle*{1}}
     \put(41,21){\circle*{1}}
     
     \put(1,1){\line(1,0){40}}
     \put(1,21){\line(1,0){40}}
     \put(1,1){\line(0,1){20}}
     \put(21,1){\line(0,1){20}}
     \put(41,1){\line(0,1){20}}
     
     \put(1,1){\line(1,1){8}}
     \put(21,1){\line(1,1){8}}
     \put(21,21){\line(-1,-1){8}}
     \put(41,21){\line(-1,-1){8}}
     
     \end{picture}
    } 
\put(55,46){
     \begin{picture}(47,16)

     \put(1,1){\circle*{1}}
     \put(1,21){\circle*{1}}
     \put(21,1){\circle*{1}}
     \put(21,21){\circle*{1}}
     \put(41,1){\circle*{1}}
     \put(41,21){\circle*{1}}
     
     \put(1,1){\line(1,0){40}}
     \put(1,21){\line(1,0){40}}
     \put(1,1){\line(0,1){20}}
     \put(21,1){\line(0,1){20}}
     \put(41,1){\line(0,1){20}}
     
     \put(1,21){\line(1,-1){20}}
     \put(21,1){\line(1,1){20}}
     \end{picture}
    } 
\put(55,9){
     \begin{picture}(47,16)

     \put(1,1){\circle*{1}}
     \put(1,21){\circle*{1}}
     \put(21,1){\circle*{1}}
     \put(21,21){\circle*{1}}
     \put(41,1){\circle*{1}}
     \put(41,21){\circle*{1}}
     
     \put(1,1){\line(1,0){40}}
     \put(1,21){\line(1,0){40}}
     \put(1,1){\line(0,1){20}}
     \put(21,1){\line(0,1){20}}
     \put(41,1){\line(0,1){20}}
     
     \put(1,21){\line(1,-1){20}}
     \put(21,21){\line(1,-1){20}}
     \end{picture}
    } 
\put(0,9){
     \begin{picture}(47,16)

     \put(1,1){\circle*{1}}
     \put(1,21){\circle*{1}}
     \put(21,1){\circle*{1}}
     \put(21,21){\circle*{1}}
     \put(41,1){\circle*{1}}
     \put(41,21){\circle*{1}}
     
     \put(1,1){\line(1,0){40}}
     \put(1,21){\line(1,0){40}}
     \put(1,1){\line(0,1){20}}
     \put(21,1){\line(0,1){20}}
     \put(41,1){\line(0,1){20}}
     
     \put(1,1){\line(1,1){20}}
     \put(21,21){\line(1,-1){20}}
     \end{picture}
    } 
\put(46,57){\vector(1,0){10}}
\put(78,43){\vector(0,-1){10}}
\put(56,19){\vector(-1,0){10}}
\put(23,31){\vector(0,1){10}}
\put(11.5,55){$e_1$}
\put(31.5,55){$e_2$}
\put(14,49){$\Delta_1$}
\put(25,60){$\Delta_2$}
\put(5,60){$\Delta_3$}
\put(34,49){$\Delta_4$}
\put(21,42){$x$}
\put(21,70){$y$}

\put(51,0){c}
\end{picture}

}
\put(93,0){Fig. 3}
\end{picture}

By $\mathcal S_3^1(p,q)$ denote the set of all homology classes
$\bar\alpha_3(L,e_1,e_2)$ such that:\\
1) the triangle $\Delta_1$ has a unique common vertex $x$ with the
triangle $\Delta_2$, see Fig.~3,~b;\\
2) there are exactly $p$ triangles containing the vertex $x$ and
situated in the angle~$\vartheta_1$;\\
3) there are exactly $q$ triangles containing the vertex $x$ and
situated in the angle~$\vartheta_2$.

By $\mathcal S_3^2(p,q)$ denote the set of all homology classes
$\bar\alpha_3(L,e_1,e_2)$ such that:\\
1) the triangles $\Delta_1$ and $\Delta_2$ have a  common
edge with endpoints $x$ and $y$, see Fig. 3, c;\\
2) there are exactly $p$ triangles that contain the vertex $x$
and coincide with none of the triangles $\Delta_1$, $\Delta_2$, and $\Delta_4$;\\
3) there are exactly $q$ triangles that contain the vertex $y$ and
coincide  with none of the triangles $\Delta_1$, $\Delta_2$, and
$\Delta_3$.

Let $x,y,z$ be vertices of $L$. Suppose there exists the vertex
$u$ such that $\{ u,x,y\}$, $\{ u,y,z\}$, and $\{ u,z,x\}$ are
triangles of $L$. Then by $\alpha_4(L,x,y,z)$ denote the cycle
shown in Fig.~4. Let $\mathcal S_4(p,q,r)$ be the set of all
homology classes $\bar\alpha_4(L,x,y,z)$ such that there are
exactly $p$, $q$, and $r$ triangles that contain respectively the
vertices $x$, $y$, and $z$ and coincide with none of the
triangles $\{ u,x,y\}$, $\{ u,y,z\}$, and $\{u,z,x\}$.

\begin{picture}(200,75)
\put(10,0){
\begin{picture}(72,75)

\put(20,49){
     \begin{picture}(32,22)
     \put(1,1){\circle*{1}}
     \put(16,8.5){\circle*{1}}
     \put(16,21){\circle*{1}}
     \put(31,1){\circle*{1}}
     \put(1,1){\line(1,0){30}}
     \put(1,1){\line(3,4){15}}
     \put(1,1){\line(2,1){15}}
     \put(16,8.5){\line(0,1){12.5}}
     \put(16,8.5){\line(2,-1){15}}
     \put(16,21){\line(3,-4){15}}
     \end{picture}
    } 
\put(40,9){
     \begin{picture}(32,28)
     \put(1,1){\circle*{1}}
     \put(11,11){\circle*{1}}
     \put(16,26){\circle*{1}}
     \put(21,11){\circle*{1}}
     \put(31,1){\circle*{1}}
     \put(1,1){\line(1,0){30}}
     \put(1,1){\line(3,5){15}}
     \put(1,1){\line(1,1){10}}
     \put(11,11){\line(1,0){10}}
     \put(11,11){\line(1,3){5}}
     \put(16,26){\line(3,-5){15}}
     \put(16,26){\line(1,-3){5}}
     \put(21,11){\line(1,-1){10}}
     \put(11,11){\line(2,-1){20}}
     \end{picture}
    } 
\put(0,9){
     \begin{picture}(32,28)
     \put(1,1){\circle*{1}}
     \put(11,11){\circle*{1}}
     \put(16,26){\circle*{1}}
     \put(21,11){\circle*{1}}
     \put(31,1){\circle*{1}}
     \put(1,1){\line(1,0){30}}
     \put(1,1){\line(3,5){15}}
     \put(1,1){\line(1,1){10}}
     \put(11,11){\line(1,0){10}}
     \put(11,11){\line(1,3){5}}
     \put(16,26){\line(3,-5){15}}
     \put(16,26){\line(1,-3){5}}
     \put(21,11){\line(1,-1){10}}
     \put(1,1){\line(2,1){20}}
     \end{picture}
    }

\put(41,20){\vector(-1,0){10}}
\put(43,45){\vector(1,-2){5}}
\put(23,35){\vector(1,2){5}}
\put(19,48){$x$}
\put(37,72){$y$}
\put(53,48){$z$}
\put(34,58){$u$}

\put(30,0){Fig. 4}
\end{picture}

} 
\put(110,0){
\begin{picture}(101,75)

\put(4,49){
     \begin{picture}(22,22)
     \put(1,1){\circle*{1}}
     \put(21,1){\circle*{1}}
     \put(1,21){\circle*{1}}
     \put(21,21){\circle*{1}}
     \put(1,1){\line(1,0){20}}
     \put(1,1){\line(0,1){20}}
     \put(1,21){\line(1,0){20}}
     \put(21,1){\line(0,1){20}}
     \put(1,1){\line(1,1){20}}
     \end{picture}
    } 
\put(44,49){
     \begin{picture}(22,22)
     \put(1,1){\circle*{1}}
     \put(21,1){\circle*{1}}
     \put(1,21){\circle*{1}}
     \put(21,21){\circle*{1}}
     \put(11,6){\circle*{1}}
     \put(1,1){\line(1,0){20}}
     \put(1,1){\line(0,1){20}}
     \put(1,21){\line(1,0){20}}
     \put(21,1){\line(0,1){20}}
     \put(1,1){\line(1,1){20}}
     \put(1,1){\line(2,1){10}}
     \put(11,6){\line(2,3){10}}
     \put(11,6){\line(2,-1){10}}
     \end{picture}
    } 
\put(79,29){
     \begin{picture}(22,22)
     \put(1,1){\circle*{1}}
     \put(21,1){\circle*{1}}
     \put(1,21){\circle*{1}}
     \put(21,21){\circle*{1}}
     \put(11,6){\circle*{1}}
     \put(1,1){\line(1,0){20}}
     \put(1,1){\line(0,1){20}}
     \put(1,21){\line(1,0){20}}
     \put(21,1){\line(0,1){20}}
     \put(1,21){\line(2,-3){10}}
     \put(1,1){\line(2,1){10}}
     \put(11,6){\line(2,3){10}}
     \put(11,6){\line(2,-1){10}}
     \end{picture}
    } 
\put(44,9){
     \begin{picture}(22,22)
     \put(1,1){\circle*{1}}
     \put(21,1){\circle*{1}}
     \put(1,21){\circle*{1}}
     \put(21,21){\circle*{1}}
     \put(11,6){\circle*{1}}
     \put(1,1){\line(1,0){20}}
     \put(1,1){\line(0,1){20}}
     \put(1,21){\line(1,0){20}}
     \put(21,1){\line(0,1){20}}
     \put(1,21){\line(1,-1){20}}
     \put(1,1){\line(2,1){10}}
     \put(1,21){\line(2,-3){10}}
     \put(11,6){\line(2,-1){10}}
     \end{picture}
    } 

\put(4,9){
     \begin{picture}(22,22)
     \put(1,1){\circle*{1}}
     \put(21,1){\circle*{1}}
     \put(1,21){\circle*{1}}
     \put(21,21){\circle*{1}}
     \put(1,1){\line(1,0){20}}
     \put(1,1){\line(0,1){20}}
     \put(1,21){\line(1,0){20}}
     \put(21,1){\line(0,1){20}}
     \put(1,21){\line(1,-1){20}}
     \end{picture}
    } 

\put(30,60){\vector(1,0){10}}
\put(68,52){\vector(2,-1){10}}
\put(40,20){\vector(-1,0){10}}
\put(15,35){\vector(0,1){10}}
\put(78,32){\vector(-2,-1){10}}
\put(3,47){$x$}
\put(3,71){$y$}
\put(27,71){$z$}
\put(27,47){$u$}

\put(45,0){Fig. 5}
\end{picture}

}
\end{picture}

Let $x,y,z,u$ be the vertices of $L$. Suppose the full subcomplex
of $L$ spanned by the set $\{ x,y,z,u\}$ consists of the triangles
$\{ x,y,z\}$ and $\{ x,z,u\}$, their edges and vertices. Then let
$\alpha_5(L,x,y,z,u)$ be the cycle shown in Fig.~5. Let $\mathcal
S_5(p,q,r,k)$ be the set of all homology classes
$\bar\alpha_5(L,x,y,z,u)$ such that there are exactly  $p$, $q$,
$r$, and $k$ triangles that contain respectively the vertices
$x$, $y$, $z$, and $u$ and coincide neither with $\{ x,y,z\}$ nor
with $\{ x,z,u\}$.

Let $x,y,z,u,v$ be the vertices of $L$. Suppose the full
subcomplex of $L$ spanned by the set $\{ x,y,z,u,v\}$ consists of
the triangles $\{ x,y,z\}$, $\{ x,z,u\}$, and $\{ x,u,v\}$, their
edges and vertices. Then let $\alpha_6(L,x,y,z,u,v)$ be the cycle
shown in Fig.~6. Let $\mathcal S_6(p,q,r,k,l)$ be the set of all
homology classes $\bar\alpha_6(L,x,y,z,u,v)$ such that there are
exactly $p$, $q$, $r$, $k$, and $l$ triangles that contain
respectively the vertices $x$, $y$, $z$, $u$, and $v$  and
coincide with none of the triangles $\{ x,y,z\}$, $\{ x,z,u\}$,
and $\{ x,u,v\}$.

\begin{picture}(101,75)
\put(60,0){
\begin{picture}(0,0)
\put(4,49){
     \begin{picture}(22,22)
     \put(6,1){\circle*{1}}
     \put(1,11){\circle*{1}}
     \put(11,21){\circle*{1}}
     \put(16,1){\circle*{1}}
     \put(21,11){\circle*{1}}
     \put(1,11){\line(1,1){10}}
     \put(1,11){\line(1,-2){5}}
     \put(6,1){\line(1,0){10}}
     \put(11,21){\line(1,-1){10}}
     \put(16,1){\line(1,2){5}}
     \put(6,1){\line(1,4){5}}
     \put(6,1){\line(3,2){15}}
     \end{picture}
    }
\put(44,49){
     \begin{picture}(22,22)
     \put(6,1){\circle*{1}}
     \put(1,11){\circle*{1}}
     \put(11,21){\circle*{1}}
     \put(16,1){\circle*{1}}
     \put(21,11){\circle*{1}}
     \put(1,11){\line(1,1){10}}
     \put(1,11){\line(1,-2){5}}
     \put(6,1){\line(1,0){10}}
     \put(11,21){\line(1,-1){10}}
     \put(16,1){\line(1,2){5}}
     \put(1,11){\line(1,0){20}}
     \put(6,1){\line(3,2){15}}
     \end{picture}
    }
\put(79,29){
     \begin{picture}(22,22)
     \put(6,1){\circle*{1}}
     \put(1,11){\circle*{1}}
     \put(11,21){\circle*{1}}
     \put(16,1){\circle*{1}}
     \put(21,11){\circle*{1}}
     \put(1,11){\line(1,1){10}}
     \put(1,11){\line(1,-2){5}}
     \put(6,1){\line(1,0){10}}
     \put(11,21){\line(1,-1){10}}
     \put(16,1){\line(1,2){5}}
     \put(1,11){\line(1,0){20}}
     \put(1,11){\line(3,-2){15}}
     \end{picture}
    }
\put(44,9){
     \begin{picture}(22,22)
     \put(6,1){\circle*{1}}
     \put(1,11){\circle*{1}}
     \put(11,21){\circle*{1}}
     \put(16,1){\circle*{1}}
     \put(21,11){\circle*{1}}
     \put(1,11){\line(1,1){10}}
     \put(1,11){\line(1,-2){5}}
     \put(6,1){\line(1,0){10}}
     \put(11,21){\line(1,-1){10}}
     \put(16,1){\line(1,2){5}}
     \put(11,21){\line(1,-4){5}}
     \put(1,11){\line(3,-2){15}}
     \end{picture}
    }

\put(4,9){
     \begin{picture}(22,22)
     \put(6,1){\circle*{1}}
     \put(1,11){\circle*{1}}
     \put(11,21){\circle*{1}}
     \put(16,1){\circle*{1}}
     \put(21,11){\circle*{1}}
     \put(1,11){\line(1,1){10}}
     \put(1,11){\line(1,-2){5}}
     \put(6,1){\line(1,0){10}}
     \put(11,21){\line(1,-1){10}}
     \put(16,1){\line(1,2){5}}
     \put(6,1){\line(1,4){5}}
     \put(11,21){\line(1,-4){5}}
     \end{picture}
    }

\put(30,60){\vector(1,0){10}}
\put(68,52){\vector(2,-1){10}}
\put(40,20){\vector(-1,0){10}}
\put(15,35){\vector(0,1){10}}
\put(78,32){\vector(-2,-1){10}}
\put(8,47){$x$}
\put(3,60){$y$}
\put(16,71){$z$}
\put(27,60){$u$}
\put(22,47){$v$}

\put(45,0){Fig. 6}
\end{picture}
}
\end{picture}

By  $\mathcal S$ denote the union of all the sets $\mathcal
S_1^0$, $\mathcal S_1^1(p,q)$, $\mathcal S_1^2(p,q)$, $\mathcal
S_2^0$, $\mathcal S_2^1(p,q)$, $\mathcal S_2^2(p,q)$, $\mathcal
S_3^0$, $\mathcal S_3^1(p,q)$, $\mathcal S_3^2(p,q)$, $\mathcal
S_4(p,q,r)$, $\mathcal S_5(p,q,r,k)$, and $\mathcal
S_6(p,q,r,k,l)$.

\begin{propos}
The set $\mathcal S$ generates the group $H_1(\Gamma_2;\mathbb
Z)$.
\end{propos}

\subsection{The formula}

By Theorem 2.1, there exists a unique generator $\phi\in
H^4(\mathcal{T}^*(\mathbb{Q}))\cong\mathbb{Q}$ such that
$\phi^{\sharp}(M)=p_1(M)$ for any manifold $M$. The following
theorem gives an explicit description of the generator $\phi$.

\begin{theorem}
Let $c_0:\mathcal S\to\mathbb Q$ be the function
given by

\begin{multline*}
\shoveright{c_0(\bar\alpha)=0,\,\,\,\bar\alpha\in\mathcal S^0_1\cup\mathcal S^0_2\cup\mathcal S^0_3}
\end{multline*}

\begin{multline*}
\shoveright{c_0(\bar\alpha)=\frac{q-p}{(p+q+2)(p+q+3)(p+q+4)},\,\,\,\bar\alpha\in \mathcal S^1_1(p,q)\cup\mathcal S^1_2(p,q)\cup\mathcal S^1_3(p,q)}
\end{multline*}

\begin{multline*}
\shoveright{c_0(\bar\alpha)=\frac{q}{(q+2)(q+3)(q+4)}-\frac{p}{(p+2)(p+3)(p+4)},\,\,\,\bar\alpha\in \mathcal S^2_1(p,q)\cup\mathcal S^2_3(p,q)}
\end{multline*}

\begin{multline*}
\shoveright{c_0(\bar\alpha)=\frac{q}{(q+2)(q+3)(q+4)}+ \frac{p}{(p+2)(p+3)(p+4)},\,\,\,\bar\alpha\in \mathcal S^2_2(p,q)}
\end{multline*}

\begin{multline*}
\shoveright{c_0(\bar\alpha)= \frac{1}{(p+2)(p+3)} -\frac{1}{(q+2)(q+3)}+ \frac{1}{(r+2)(r+3)}-\frac1{12},\,\,\,\bar\alpha\in\mathcal S_4(p,q,r)}
\end{multline*}

\begin{multline*}
c_0(\bar\alpha)= \frac{1}{(p+2)(p+3)} -\frac{1}{(q+2)(q+3)}- \frac{1}{(r+2)(r+3)}\\
+\frac1{(k+2)(k+3)},\,\,\,\bar\alpha\in\mathcal S_5(p,q,r,k)
\end{multline*}

\begin{multline*}
c_0(\bar\alpha)= \frac{1}{(p+2)(p+3)} +\frac{1}{(q+2)(q+3)}+ \frac{1}{(r+2)(r+3)}+\frac1{(k+2)(k+3)}\\
+\frac1{(l+2)(l+3)}-\frac1{12},\,\,\,\bar\alpha\in\mathcal S_6(p,q,r,k,l)
\end{multline*}

{\noindent There exists a unique linear extension of $c_0$ to
$H_1(\Gamma_2;\mathbb Z)$. This extension, which we also  denote
by $c_0$, belongs to the group $H_{\mathbb
Z_2}^1(\Gamma_2;\mathbb Q) =\Hom_{\mathbb
Z_2}(H_1(\Gamma_2;\mathbb Z),\mathbb Q)$. Then $c_0=s^*(\phi)$.
Thus $s$ maps isomorphically the affine space of all local
formulae for the first Pontrjagin class to the affine space of
all cocycles $\widehat c_0\in C^1_{\mathbb Z_2}(\Gamma_2;\mathbb Q)$
representing the class $c_0$.}
\end{theorem}

\begin{note}
Notice that similar numbers appear in~\cite{K} in a quite different
problem. In this paper M. \`E Kazarian obtained a formula for a
Chern-Euler class of a  circle bundle over a closed surface in
terms of singularities of a generic function's restrictions to
the fibers.
\end{note}

How to make calculations using this theorem? First we need to
choose a representative $\widehat c_0$ of the class $c_0$ to fix a
local formula $f$ for the first Pontrjagin class.  Let us now
show how to calculate the value $f(L)$, where $L$ is an oriented
$3$-dimensional PL sphere. Let $\beta_1,\beta_2,\ldots,\beta_l$ be
the sequence of bistellar moves transforming the boundary of
$4$-simplex to the PL sphere $L$. Let $L_j$ be the PL sphere
obtained from $\partial\Delta^4$ by applying the bistellar moves
$\beta_1,\beta_2,\ldots,\beta_{j-1}$. Let $W_j$ be the set of all
vertices $v\in L_j$ such that the bistellar move $\beta_j$
induces an essential bistellar move $\beta_{jv}$ of $\Lk v$. Then
$$f(L)=\sum_{j=1}^l\sum_{v\in W_j}(\widehat c_0)(e_{\beta_{jv}})$$

Assume that we need to calculate the first Pontrjagin number of
an oriented $4$-dimensional combinatorial manifold $K$. Then we
don't need to choose a representative $\widehat c_0$ of the
cohomology class $c_0$. Let $L^1,L^2,\ldots,L^k$ be the links of
all vertices of $K$. Define the bistellar moves $\beta_j^i$, the
PL spheres $L^i_j$, the sets $W_j^i$ and the bistellar moves
$\beta^i_{jv}$ as above. 

\begin{corr}
The first Pontrjagin number of the
manifold $K$ is equal to
$$
\frac{1}{2}c_0\left(\sum_{i=1}^k\sum_{j=1}^{l_i}\sum_{v\in
W_j^i}(e_{\beta^i_{jv}}-\widetilde e_{\beta^i_{jv}})\right)
$$
where $\widetilde e$ is the edge of the graph $\Gamma_2$ such that
the action of the group $\mathbb Z_2$ takes $e$ to $\widetilde e$.
\end{corr}

\begin{proof}[Sketch of the proof of Theorem 3.1]
Since
$\star:H^4(\mathcal{T}^*(\mathbb{Q}))\to\Hom(\Omega_4,\mathbb Q)$
is an epimorphism, we see that $\dim H^4(\mathcal T^*(\mathbb Q))\ge 1$. On
the other hand, in section 3.4 we prove that for any cohomology
class $c\in\widetilde A^3$ there exists $\lambda\in \mathbb Q$ such
that $c(\bar\alpha)=\lambda c_0(\bar\alpha)$ for any
$\bar\alpha\in\mathcal S$. Therefore $\dim\widetilde A^3\le 1$. But
$s^*:H^4(\mathcal T^*(\mathbb Q))\to\widetilde A^3$ is a monomorphism. Hence
$\dim H^4(\mathcal T^*(\mathbb Q))=\dim\widetilde A^3=1$. Therefore the
cohomology class $c_0$ is well-defined, $c_0\in\widetilde A^3$ and
$s^*(\phi)=\lambda c_0$ for some rational $\lambda\ne 0$. In
section 3.5 we prove that $\lambda=1$.
\end{proof}

\subsection{The group $\widetilde A^3$}

Consider an arbitrary $c\in\widetilde
A^3\subset\Hom_{\mathbb{Z}_2}(H_1(\Gamma_2;\mathbb{Z}),\mathbb{Q})$.

Let $\alpha=\alpha_i(L,\ldots)$ be a cycle shown in one of
Fig.~1--6. By $X(\alpha)$ denote the set of all vertices that are
denoted in the corresponding figure by one of the letters $x$,
$y$, $z$, $u$, and $v$. The generator $\bar\alpha\in\mathcal S$
is called {\it regular} if it satisfies the following
conditions:\\
1) $\bigcup\limits_{a\in X(\alpha)}\St a$ is a full subcomplex of
$L$.\\
2) If $w\notin X(\alpha)$, $a,b\in X(\alpha)$, and
$\{w,a\},\{w,b\}\in L$, then $\{ w,a,b\}\in L$.

Consider $\bar\alpha_1(L,\Delta_1,\Delta_2)\in\mathcal S^0_1$. Let
$K\in\mathcal T_4$ be a PL sphere containing a vertex $u$ such that $\Lk
u\cong L$ and $\St u$ is a full subcomplex of $K$. Identify the
two simplicial complexes  $\Lk u$ and $L$. Put
$\widetilde\Delta_1=\Delta_1\cup\{ u\}$ and
$\widetilde\Delta_2=\Delta_2\cup\{ u\}$. Apply to $K$ the bistellar
move associated with $\widetilde\Delta_1$; by $z_1$ denote the new
vertex created by this bistellar move. Then apply to the obtained
PL sphere the bistellar move associated with $\widetilde\Delta_2$; by
$z_2$ denote the new vertex created by this bistellar move.
Apply  to the obtained PL sphere the bistellar move associated
with the vertex $z_1$. Finally, apply to the obtained PL sphere
the bistellar move associated with the vertex $z_2$. By $\gamma$
denote the obtained cycle in the graph $\Gamma_3$. For any vertex
$v$ of $K$ the sequence of bistellar moves transforming $K$ to
itself induces the sequence of bistellar moves transforming $\Lk
v$ to itself. Hence the cycle $\gamma$ induces the cycle
$\beta_v$ in the graph $\Gamma_2$. Then
$\delta^*(c)(\bar\gamma)=\sum_vc(\bar\gamma_v)$. Now
$\delta^*(c)=0$ because $c\in\widetilde A^3$. For any vertex $v\in K$
distinct from the vertex $u$ the cycle $\gamma_v$ is homologous
to zero. The cycle $\gamma_u$ coincides with the cycle
$\alpha_1(L,\Delta_1,\Delta_2)$. Hence,
$c(\bar\alpha_1(L,\Delta_1,\Delta_2))=0$.

Similarly, $c(\bar\alpha)=0$ for any $\bar\alpha\in\mathcal S_2^0\cup\mathcal S_3^0$.

\begin{propos}
For any $p,q>0$ the function $c$ is a constant function on the
set $\mathcal{S}_1^1(p,q)$.
\end{propos}

\begin{proof}
Consider $\bar\alpha_1(L^{(i)},\Delta^{(i)}_1,\Delta^{(i)}_2)\in
\mathcal S_1^1(p,q)$, $i=1,2$. Let $x^{(i)}$ be the common vertex
of the triangles $\Delta^{(i)}_1$ and $\Delta^{(i)}_2$. Assume
that $\bar\alpha_1(L^{(1)},\Delta^{(1)}_1,\Delta^{(1)}_2)$ is
regular. Obviously, there exists a $3$-dimensional oriented PL
sphere $K$ containing an edge $e$ with
endpoints $u$ and $w$ such that the following conditions hold:\\
1) The link of $e$ is a $(p+q+2)$-gon containing two  edges denoted by $e_1$ and $e_2$.\\
2) The link of $u$ is isomorphic to $L^{(1)}$. This isomorphism
takes the triangle spanned by the vertex $w$ and the edge $e_j$
to the triangle $\Delta^{(1)}_j$, $j=1,2$.\\
3) The link of $w$ is isomorphic to $-L^{(2)}$. This isomorphism
takes the triangle spanned by the vertex $u$ and the edge $e_j$
to the triangle $\Delta^{(2)}_j$, $j=1,2$.

By $\widetilde\Delta_j$ denote the $3$-simplex of $L$ spanned by the
edges $e$ and $e_j$. Let the cycle $\gamma$ be as above. For any
vertex $v\in L$ distinct from the vertices $u$ and $w$ the induced
cycle $\gamma_v$ is homologous to zero. The cycle $\gamma_u$
coincides with the cycle
$\alpha_1(L^{(1)},\Delta^{(1)}_1,\Delta^{(2)}_2)$, the cycle
$\gamma_v$ coincides with the cycle
$\alpha_1(-L^{(2)},\Delta_1^{(2)},\Delta_2^{(2)})$. Hence
$c(\bar\alpha_1(L^{(1)},\Delta^{(1)}_1,\Delta^{(2)}_2))=
c(\bar\alpha_1(L^{(2)},\Delta_1^{(2)},\Delta_2^{(2)}))$. To
conclude the proof note that each set $\mathcal{S}_1^1(p,q)$
contains a regular generator.
\end{proof}

By  $\rho(p,q)$ denote the value of the function $c$ on the set
$\mathcal{S}_1^1(p,q)$.

\begin{propos}
For any $p,q>0$ the function $c$ is a constant function on the
set $\mathcal{S}_1^2(p,q)$. Let $\tau(p,q)$ be the value of the
function $c$ on the set $\mathcal{S}_1^2(p,q)$. Then $\tau
(p,q)+\tau (q,r)+\tau (r,p)=0$ for any $p,q,r>0$.
\end{propos}

\begin{proof}
Consider
$\bar\alpha_i(L^{(1)},\Delta^{(1)}_1,\Delta^{(1)}_2)\in\mathcal
S_1^2(p,q)$. Consider $r>0$. Choose regular generators
$\bar\alpha_1(L^{(2)},\Delta_1^{(2)},\Delta_2^{(2)})\in\mathcal
S_1^2(q,r)$ and
$\bar\alpha_1(L^{(3)},\Delta_1^{(3)},\Delta_2^{(3)})\in\mathcal
S_1^2(r,p)$. There exists a $3$-dimensional oriented PL sphere $K$
containing a $2$-simplex $\Delta_0$ with vertices $u^{(1)}$,
$u^{(2)}$, and $u^{(3)}$ such that the following conditions hold:\\
1) The link of $u^{(i)}$ is isomorphic to $L^{(i)}$.\\
2) This isomorphism takes the $2$-dimensional face of the simplex
$\widetilde\Delta_j$ opposite the vertex $u^{(i)}$ to the triangle
$\Delta_j^{(i)}$, $j=1,2$, where
$\widetilde\Delta_1,\widetilde\Delta_2\in L$ are the two tetrahedrons
containing $\Delta_0$.

As in the prove of Proposition 3.5 we obtain
$$c(\bar\alpha_1(L^{(1)},\Delta^{(1)}_1,\Delta^{(1)}_2))+
c(\bar\alpha_1(L^{(2)},\Delta^{(2)}_1,\Delta_2^{(2)}))+
c(\bar\alpha_1(L^{(3)},\Delta_1^{(3)},\Delta_2^{(3)}))=0$$
Instead of $\bar\alpha_1(L^{(1)},\Delta^{(1)}_1,\Delta^{(1)}_2)$ we can take any
generator $\bar\alpha\in\mathcal S_1^2(p,q)$. Hence the function $c$
is constant on $\mathcal{S}_1^2(p,q)$. The equality $\tau
(p,q)+\tau (q,r)+\tau (r,p)=0$ follows.
\end{proof}

Evidently, $\tau(p,q)=-\tau(q,p)$. Therefore there exists a
function $\chi : \mathbb{Z}_{>0}\to\mathbb{Q}$ such that $\tau
(p,q)=\chi (q)-\chi (p)$ for any $p,q>0$. The function $\chi$ is
unique up to a constant. Extend the function
$\rho$ to the function of $\mathbb Z_{\ge 0}\times\mathbb Z_{\ge
0}$ to $\mathbb Q$ by putting $\rho(0,p)=\chi(p)$,
$\rho(p,0)=-\chi (p)$, and $\rho (0,0)=0$.

\begin{propos}
The function $\rho$ satisfies the following equations.
\begin{multline*}
(i)\;\;\rho(p,q)=-\rho (q,p)\\ \shoveleft{(ii)\;\rho(p,q+r+2)+\rho(q,r+p+2)+\rho(r,p+q+2)=\rho(p,q+r+1)}\\ +\rho(q,r+p+1)+\rho(r,p+q+1)
\end{multline*}
\end{propos}

\begin{proof}
Equation $(i)$ follows from
$\alpha_1(L,\Delta_1,\Delta_2)=-\alpha_1(L,\Delta_2,\Delta_1)$.

Let $L$ be an oriented simplicial $2$-sphere containing a vertex
$x$ such that there exist exactly $(p+q+r+3)$ triangles of $L$
containing $x$. Let us go round the vertex $x$ clockwise. Let
$\Delta_1$, $\Delta_2$, and $\Delta_3$ be triangles containing
the vertex $x$ such that we pass through the triangle $\Delta_1$,
then through $r$ other triangles, then through the triangle
$\Delta_2$, then through $p$ other triangles, then through the
triangle  $\Delta_3$, then through $q$ other triangles, and then
again through the triangle $\Delta_1$. By $L_j$ denote the PL
sphere obtained from $L$ by applying the bistellar move associated
with $\Delta_j$. It is easy to prove that
$\sum_{j=1}^3\alpha_1(L_j,\Delta_{j+1},\Delta_{j+2})=
\sum_{j=1}^3\alpha_1(L,\Delta_{j+1},\Delta_{j+2})$, where the sums
of indices are modulo $3$. Applying $c$ to the homology classes
of both parts of this equality, we obtain equation $(ii)$.
\end{proof}

\begin{propos}Suppose the function $\rho : \mathbb{Z}_{\ge
0}\times\mathbb{Z}_{\ge 0}\to\mathbb{Q}$ satisfies equations
$(i)$ and $(ii)$. Then there exist $b_1\in\mathbb Q$ and
$\lambda\in\mathbb Q$ such that $\rho(p,q)=\frac{\lambda
(q-p)}{(p+q+2)(p+q+3)(p+q+4)}$ for any $p,q>0$ and
$\rho(0,q)=\frac{\lambda q}{(q+2)(q+3)(q+4)}+b_1$ for any $q>0$.
\end{propos}

Without loss of generality we may assume that
$\rho(p,q)=\frac{\lambda (q-p)}{(p+q+2) (p+q+3)(p+q+4)}$ for any
$p,q\ge 0$.

The proofs of the following propositions are similar two the
proofs of Propositions 3.5 and 3.6.

\begin{propos}
$c(\bar\alpha)=\rho(p,q)$ for any $\bar\alpha\in\mathcal
S_2^1(p,q)\cup\mathcal S_3^1(p,q)$.
\end{propos}
\begin{propos}
There exists a constant $b_2\in\mathbb Q$ such that
$c(\bar\alpha)=\rho(0,q)+\rho(0,p)+b_2$ for any
$\bar\alpha\in\mathcal S_2^2(p,q)$.
\end{propos}
\begin{propos}
$c(\bar\alpha)=\rho(0,q)-\rho(0,p)$ for any $\bar\alpha\in\mathcal
S_3^2(p,q)$.
\end{propos}

Let $p_{12}$, $p_{13}$, $p_{14}$, $p_{23}$, $p_{24}$, $p_{25}$, $p_{34}$, $p_{35}$, $p_{36}$,
$p_{45}$, $p_{46}$ и $p_{56}$ be integers greater than $2$. Let
$\omega_j\in Z_1(\Gamma_2;\mathbb{Z})$, $j=1,2,\ldots,6$ be cycles such that the following
conditions hold:\\
1) $\omega_1=\alpha_4(L_1,x_1,y_1,z_1)$, $\bar\omega_1\in\mathcal S_4(p_{13},p_{14},p_{12})$.\\
2) $\omega_2=\alpha_5(L_2,x_2,y_2,z_2,u_2)$, $\bar\omega_2\in\mathcal S_5(p_{23},p_{12},p_{24},p_{25})$.\\
3) $\omega_3=\alpha_6(L_3,x_3,y_3,z_3,u_3,v_3)$, $\bar\omega_3\in\mathcal S_6(p_{34},p_{13},p_{23},p_{35},p_{36})$.\\
4) $\omega_4=\alpha_6(L_4,x_4,y_4,z_4,u_4,v_4)$, $\bar\omega_4\in\mathcal S_6(p_{34},p_{46},p_{45},p_{24},p_{14})$.\\
5) $\omega_5=\alpha_5(L_5,x_5,y_5,z_5,u_5)$, $\bar\omega_5\in\mathcal S_5(p_{45},p_{56},p_{35},p_{25})$.\\
6) $\omega_6=\alpha_4(L_6,x_6,y_6,z_6)$, $\bar\omega\in\mathcal
S_4(p_{46},p_{36},p_{56})$;\\
7) At least $5$ of the generators $\bar\omega_j$ are regular.

By $L_j^{(1)}$, $j=1,6$ denote the PL sphere obtained from $L_j$
by the bistellar move associated with the vertex $u_j$.

There exists a $3$-dimensional oriented PL sphere $K$ containing
vertices $w_j$, $j=1,2,\ldots 6$ such that the following conditions hold:\\
1) The full subcomplex of $K$ spanned by the set $\{
w_j,\,j=1,\ldots 6\}$ consists of the tetrahedrons $\{
w_1,w_2,w_3,w_4\}$, $\{ w_2,w_3,w_4,w_5\}$, and
$\{ w_3,w_4,w_5,w_6\}$ and all their faces.\\
2) For any edge $\{ w_i,w_j\}\in L$, $i<j$  there are exactly
$p_{ij}$ tetrahedrons that contain the edge $\{ w_i,w_j\}$ and
coincide with none of the tetrahedrons $\{ w_1,w_2,w_3,w_4\}$,
$\{ w_2,w_3,w_4,w_5\}$, and $\{ w_3,w_4,w_5,w_6\}$.\\
3) The links of the vertices $w_j$ can be identified with the
complexes $L_j$ for $1<j<6$ and $L_j^{(1)}$ for $j=1,6$ so that
the vertices will be identified in the following way:
$w_1=y_2=y_3=v_4$, $w_2=z_1=z_3=u_4=u_5$,
$w_3=x_1=x_2=x_4=z_5=y_6$, $w_4=y_1=z_2=x_3=x_5=x_6$,
$w_5=u_2=u_3=z_4=z_6$, $w_6=v_3=y_4=y_5$.

Apply to $L$ the following sequence of bistellar moves:\\
1) Replace the $3$-simplices $\{ w_1,w_2,w_3,w_4\}$ and $\{ w_2,w_3,w_4,w_5\}$ with the
$3$-simplices $\{ w_1,w_2,w_3,w_5\}$, $\{ w_1,w_3,w_4,w_5\}$, and $\{ w_1,w_4,w_2,w_5\}$.\\
2) Replace the $3$-simplices $\{ w_1,w_3,w_4,w_5\}$ and $\{ w_3,w_4,w_5,w_6\}$ with the
$3$-simplices $\{ w_1,w_3,w_4,w_6\}$, $\{ w_1,w_4,w_5,w_6\}$, and $\{ w_1,w_5,w_3,w_6\}$.\\
3) Replace the $3$-simplices $\{ w_1,w_2,w_3,w_5\}$ and $\{ w_1,w_5,w_3,w_6\}$ with the
$3$-simplices $\{ w_2,w_3,w_1,w_6\}$, $\{ w_2,w_1,w_5,w_6\}$, and $\{ w_2,w_5,w_3,w_6\}$.\\
4) Replace the $3$-simplices $\{ w_1,w_4,w_2,w_5\}$, $\{ w_1,w_6,w_4,w_5\}$, and
$\{ w_1,w_2,w_6,w_5\}$ with the $3$-simplices $\{ w_1,w_2,w_6,w_4\}$ and $\{ w_2,w_6,w_4,w_5\}$.\\
5) Replace the $3$-simplices $\{ w_1,w_2,w_3,w_6\}$, $\{ w_1,w_3,w_4,w_6\}$, and
$\{ w_1,w_4,w_2,w_6\}$ with the $3$-simplices $\{ w_1,w_2,w_3,w_4\}$ and $\{ w_2,w_3,w_4,w_6\}$.\\
6) Replace the $3$-simplices $\{ w_2,w_3,w_4,w_6\}$, $\{ w_2,w_4,w_5,w_6\}$, and
$\{ w_2,w_5,w_3,w_6\}$ with the $3$-simplices $\{ w_2,w_3,w_4,w_5\}$ and $\{ w_3,w_4,w_5,w_6\}$.

This sequence of bistellar moves transforms $L$ to itself. Let
$\gamma$ be the obtained cycle in the graph $\Gamma_3$. The cycle
$\gamma_{w_j}$ is homologous to the cycle  $\omega_j$ for
$j=3,5,6$. The cycle $\gamma_{w_j}$ is homologous to the cycle
$-\omega_j$ for $j=1,2,4$. Therefore,
$$c(\bar\omega_1)+c(\bar\omega_2)-c(\bar\omega_3)+c(\bar\omega_4)-c(\bar\omega_5)-c(\bar\omega_6)=0$$

Consequently the function $c$ is constant on each of the sets
$\mathcal S_4(p,q,r)$, $\mathcal S_5(p,q,r,k)$, and $\mathcal
S_6(p,q,r,k,l)$, $p,q,r,k,l\ge 3$. By $\eta(p,q,r)$,
$\zeta(p,q,r,k)$, and $\theta(p,q,r,k,l)$ respectively denote the
values of the function $c$ on these sets. {\sloppy

}

\begin{propos} There exists a constant $b_5\in\mathbb Q$ such that for any
$p,q,r,k,l\ge 3$ we have
\begin{multline*}
\theta(p,q,r,k,l)= \frac{\lambda}{(p+2)(p+3)} +\frac{\lambda}{(q+2)(q+3)}+ \frac{\lambda}{(r+2)(r+3)}\\ + \frac{\lambda}{(k+2)(k+3)}+ \frac{\lambda}{(l+2)(l+3)}+b_5
\end{multline*}
\end{propos}

\begin{proof}
Consider $\bar\alpha_6(L,x,y,z,u,v)\in\mathcal S_6(p,q,r,k,l)$,
$p,q,r,k,l\ge 3$. Let $\Delta\in L$ be a triangle such that $x$ is
a vertex of $\Delta$, $y$, $z$, $u$, and $v$ are not vertices of
$\Delta$. Let us go round the vertex $x$ clockwise. Suppose we
pass through the triangle $\Delta$, then through $p'$ other
triangles, then through the triangles  $\{ x,y,z\}$, $\{
x,z,u\}$, $\{ x,u,v\}$, then through $p''$ other triangles, and
then again through the triangle $\Delta$. Then $p'+p''=p-1$. The
cycle $\alpha_6(L,x,y,z,u,v)$ is the sequence of $5$ bistellar
moves. By $L_j$, $j=0,1,2,3,4$  denote the PL sphere obtained from
$L$ by the first $j$ bistellar moves of this sequence. In
particular, $L_0=L$. By $L_j^{(1)}$ denote the simplicial complex
obtained from $L$ by the bistellar move associated with the
triangle $\Delta$. Let us define the graph $G$ in the following
way. The vertex set of $G$ is the set $\{
L_0,\ldots,L_4,L_0^{(1)},\ldots,L_4^{(1)}\}$. To the $j$-th
bistellar move of the cycle $\alpha_6(L,x,y,z,u,v)$ assign the
edge with endpoints $L_{j-1}$ and $L_j$ (or $L_4$ and $L_0$ if
$j=5$). Similarly, to the $j$-th bistellar move of the cycle
$\alpha_6(L_0^{(1)},x,y,z,u,v)$ assign the edge with endpoints
$L_{j-1}^{(1)}$ and $L_j^{(1)}$ (or $L_4^{(1)}$ and $L_0^{(1)}$ if
$j=5$). For any $j=0,1,2,3,4$ consider the bistellar move
associated with the triangle $\Delta$. To this bistellar move
assign the edge with the endpoints $L_j$ and $L_j^{(1)}$. There is
the canonical map of the graph $G$ to the graph $\Gamma_2$. (This
map is not necessarily an injection.) The graph $G$ is isomorphic
to a $1$-skeleton of a pentagonal prism. Let us go the round of some
$2$-dimensional face of this sphere. We obtain the cycle in the
graph $\Gamma_2$. In this way, we obtain the cycles
$\alpha_6(L,x,y,z,u,v)$, $-\alpha_6(L_0^{(1)},x,y,z,u,v)$,
$\alpha_2(L_0,\Delta,\{ x,z\} )$, $\alpha_2(L_1,\Delta,\{ u,x\}
)$,   $\alpha_2(L_2,\Delta,\{ y,u\} )$, $\alpha_2(L_3,\Delta,\{
v,y\} )$,   $\alpha_2(L_4,\Delta,\{ z,v\} )$. The sum of all
these cycles is equal to zero. Consequently,
$$\theta(p,q,r,k,l)-\theta(p+1,q,r,k,l)-\rho(p',p''+1)+\rho(p'+1,p'')=0$$
Hence,
$$\theta(p+1,q,r,k,l)-\theta(p,q,r,k,l)=-\frac{2\lambda}{(p+2)(p+3)(p+4)}$$

In addition, the function $\theta$ is cyclically symmetric. The proposition follows from the
equality
$$\sum_{i=1}^j\frac{1}{i(i+1)(i+2)}=\frac{1}{4}-\frac{1}{2(j+1)(j+2)}$$
\end{proof}

\begin{propos}
$b_2=0$.
\end{propos}

\begin{proof} Consider $\bar\alpha_6(L,x,y,z,u,v)\in\mathcal
S_6(p,q,r,k,l)$, $p,q,r,k,l\ge 3$. Let $\Delta$ be the triangle
that contains the edge $\{ x,y\}$ and does not coincide with the
triangle $\{ x,y,z\}$. Arguing as in the proof of Proposition
3.12, we see that

\begin{multline*}
\theta(p+1,q+1,r,k,l)-\theta(p,q,r,k,l)=-\frac{2\lambda}{(p+2)(p+3)(p+4)}\\-\frac{2\lambda}{(q+2)(q+3)(q+4)}-b_2
\end{multline*}

Hence, $b_2=0$.
\end{proof}

The following two propositions are proved in the same way as
Proposition~3.12.

\begin{propos}
There exists $b_3\in\mathbb Q$ such that for any $p,q,r\ge 3$ we
have
$$\eta(p,q,r)= \frac{\lambda}{(p+2)(p+3)} -\frac{\lambda}{(q+2)(q+3)}+ \frac{\lambda}{(r+2)(r+3)}+b_3$$
\end{propos}

\begin{propos} There exists $b_4\in\mathbb Q$
such that for any $p,q,r,k\ge 3$ we have

\begin{multline*}
\zeta(p,q,r,k)= \frac{\lambda}{(p+2)(p+3)} -\frac{\lambda}{(q+2)(q+3)}- \frac{\lambda}{(r+2)(r+3)}\\+ \frac{\lambda}{(k+2)(k+3)}+b_4
\end{multline*}
\end{propos}

Consider $\bar\alpha_6(L,x,y,z,u,v)\in\mathcal S_6(p,q,r,k,l)$,
$p,q,r,k,l\ge 2$. Let $\Delta$ be the triangle that contains the
edge $\{ x,y\}$ and does not coincide with the triangle $\{
x,y,z\}$. Let $L^{(1)}$ be the PL sphere obtained from $L$ by the
bistellar move associated with $\Delta$. Arguing as in the proof
of Proposition 3.13, we see that
\begin{multline*}
c(\bar\alpha_6(L,x,y,z,u,v))-c(\bar\alpha_6(L^{(1)},x,y,z,u,v))=\frac{2\lambda}{(p+2)(p+3)(p+4)}\\+\frac{2\lambda}{(q+2)(q+3)(q+4)}
\end{multline*}

Consequently for any $p,q,r,k,l\ge 2$ the function $c$ is
constant on $S_6(p,q,r,k,l)$ and the formula of Proposition 3.12
holds. Similarly the formulae of Propositions 3.14 and 3.15 hold
for any $p,q,r\ge 2$ and $p,q,r,k\ge 2$ respectively. The formula
of Proposition 3.14 holds also for $p=q=r=1$. {\sloppy

}

Obviously, $\eta(1,1,1)=0$ and $\zeta(2,2,2,2)=0$. Hence
$b_3=-\frac{\lambda}{12}$ and $b_4=0$. It is easy to prove that
$\theta(2,2,2,2,2)=-5\eta(2,2,2)=\frac{\lambda}{6}$. Hence,
$b_5=-\frac{\lambda}{12}$.

Consequently $c(\bar\alpha)=\lambda c_0(\bar\alpha)$ for any
generator $\bar\alpha\in\mathcal S$.

\subsection{The constant $\lambda$}

We have $s^*(\phi)=\lambda c_0$ for some $\lambda\in\mathbb Q$.
To prove that $\lambda=1$ we need to check that the formula of
the Corollary 3.2 holds for some oriented $4$-dimensional
combinatorial manifold with nonzero first Pontrjagin number.

In~\cite{KB} W. K\"uhnel and T. F. Banchoff constructed the
triangulation $\mathbb{C}P^2$ with $9$ vertices (see
also~\cite{KL}). The links of all vertices of this triangulation
are isomorphic to the same $3$-dimensional oriented PL sphere
$L$. ($L$ is one of the two $3$-dimensional PL spheres that are
not polytopal spheres, see~\cite{GS}). One can number the vertices
of $L$ so that the $3$ dimensional simplices are given by:
$$
\begin{array}{ccccccccccccc}
1243&&&&3476&&&&5386&&&&7165\\ 1237&&&&3465&&&&4285&&&&1785\\
1276&&&&4576&&&&4875&&&&1586\\
2354&&&&2385&&&&4817&&&&1682\\
2376&&&&2368&&&&4371&&&&1284
\end{array}
$$

Consider the following sequence of $9$ bistellar moves. This sequence transforms
$L$ to the boundary of $4$-dimensional simplex.\\
1) Replace the simplices $1243$, $1237$, and $4371$ with the simplices $1247$ and $3274$.\\
2) Replace the simplices $2354$, $2385$, and $4285$ with the simplices $2384$ and $3584$.\\
3) Replace the simplices $7165$, $1785$, and $1586$ with the simplices $1786$ and $5687$.\\
4) Replace the simplices $1786$, $1682$, and $1276$ with the simplices $1278$ and $6287$.\\
5) Replace the simplices $1247$, $1278$, $1284$, and $4817$ with the simplex $2487$.\\
6) Replace the simplices $3274$, $2384$, and $2487$ with the simplices $2387$ and $3487$.\\
7) Replace the simplices $2387$, $6287$, $2376$, and $2368$ with the simplex $6387$.\\
8) Replace the simplices $3465$, $5386$, and $3584$ with the simplices $4386$ and $5486$.\\
9) Replace the simplices $4386$, $6387$, $3476$, and $3487$ with the simplex $6487$.\\

Using this sequence of bistellar moves the formula of Corollary
3.2 can be checked by a direct calculation.

\subsection{Denominators of local 
formulae's
values}

Denote by $\mathcal{T}_{n,l}$ the set of all oriented
$(n-1)$-dimensional PL spheres with not more than $l$ vertices.
Let $f\in\mathcal T^n(\mathbb Q)$ be a local formula. By $\den_l(f)$ denote
the least common multiple of the denominators of the values
$f(L)$, where $L$ runs over all elements of $\mathcal{T}_{n,l}$.

\begin{propos}
Let $\psi\in H^n(\mathcal T^*(\mathbb Q))$ be an arbitrary
cohomology class. Then there exists a local formula $f$ representing
the class $\psi$ and the integer constant $C\ne 0$ such that the
number $\den_l(f)$ is a divisor of $C(l+1)!$ for any $l\ge n$.
\end{propos}
\begin{propos}
Let $f$ be an arbitrary local formula representing the generator
$\phi$ of the group $H^4(\mathcal{T}^*(\mathbb Q))$. Then the
number $\den_l(f)$ is divisible by the least common multiple of
the numbers $1,2,3,\ldots,l-3$ for any even $l\ge 10$.
\end{propos}
The proof of Proposition 3.16 uses Theorem 2.1 and results
of~\cite{LR}. Proposition 3.17 is a corollary of Theorem 3.1.
\begin{corr}
$H^4(\mathcal{T}^*(G))=0$ for any subgroup
$G\varsubsetneq\mathbb{Q}$.
\end{corr}

I am grateful to V. M. Buchstaber for statements of the problems and for 
constant attention to my work. I am grateful to L. A. Alania, I. V. Baskakov, 
I. A. Dynnikov, M. \`E. Kazarian, and A. S. Mischenko for usefull discussions.

 \end{document}